\documentclass{rspublic}
\usepackage{amsmath,amssymb,amsthm,color,epsfig}
\usepackage{amsbsy}
\renewcommand{\thesection}{{\arabic{section}}}

\newtheorem{Theorem}{Theorem}[section]
\newtheorem{Lemma}{Lemma}[section]

\numberwithin{equation}{section}
\newcommand{\eps}{\epsilon}
\newcommand{\mueff}{\mu_\text{eff}}
\newcommand{\epseff}{\eps_\text{eff}}
\newcommand{\omegap}{\omega_\text{\tiny p}}
\newcommand{\epsp}{{\eps_\text{\tiny p}}}
\newcommand{\mup}{{\mu_\text{\tiny p}}}
\newcommand{\xx}{{\mathbf x}}
\newcommand{\yy}{{\mathbf y}}

\newcommand{\kk}{{\hat \kappa}}

\newcommand{\epspbar}{\eps_{\bar{\text{\tiny p}}}}
\newcommand{\mupbar}{\mu_{\bar{\text{\tiny p}}}}
\newcommand{\pbar}{\bar{\text{p}}}
\newcommand{\pp}{\text{p}}

\newcommand{\RR}{\mathbb{R}}

\newcommand{\QQ}{{Q}}
\newcommand{\at}[1]{|_{\scriptscriptstyle{#1}}}
\newcommand{\nn}{{\mathbf n}}
\newcommand{\grady}{\nabla\hspace{-2pt}_\yy}

\newcommand{\kdotgrad}{\kk\hspace{-1pt}\cdot\hspace{-3pt}\nabla}

\newcommand{\dP}{{\partial P}}
\newcommand{\Pbar}{{\bar P}}
\newcommand{\honepnorm}[1]{\Vert {#1} \Vert_{H_1(P)}}
\newcommand{\honepbarnorm}[1]{\Vert #1 \Vert_{H_1(\bar P)}}
\newcommand{\ltwopnorm}[1]{\Vert #1 \Vert_{L_2(P)}}

\newcommand{\integraloverp}[1]{\langle #1 \rangle_{P}}
\newcommand{\integraloverpbar}[1]{\langle #1 \rangle_{\bar P}}
\newcommand{\integraloverq}[1]{\langle #1 \rangle_{\QQ}}
\newcommand{\vectorintegraloverp}[2]{\langle #1\! \cdot\! #2\rangle_{P}}
\newcommand{\vectorintegraloverpbar}[2]{\langle #1\! \cdot\! #2\rangle_{\bar P}}

\renewcommand{\Re}{\text{Re}\,}

\newcommand{\Omegap}{\mathbf{\Omega_{\Pbar}}}

\newcommand{\xizero}{|\xi^2_0|}

\newcommand{\thetap}{\sqrt{\theta_P}}
\newcommand{\thetapbar}{\sqrt{\theta_{\Pbar}}}

\begin{document}
\title[Sub-Wavelength Plasmonic Crystals]{Sub-Wavelength Plasmonic Crystals: \\Dispersion Relations and Effective Properties}
\author[Fortes, Lipton \& Shipman]{Santiago P. Fortes, Robert P. Lipton and Stephen P. Shipman}
\affiliation{Louisiana State University, Baton Rouge, LA, USA}
\label{firstpage}
\maketitle



\begin{abstract}{Meta-material, Plasmonic crystal, Dispersion relation, Effective property, Series solution, Catalan number}
We obtain a convergent power series expansion for the first branch of the dispersion relation for subwavelength plasmonic crystals consisting of plasmonic rods with frequency-dependent dielectric permittivity embedded in a host medium with unit permittivity. The expansion parameter is $\eta=kd=2\pi d/\lambda$, where $k$ is the norm of a fixed wavevector, $d$ is the period of the crystal and $\lambda$ is the wavelength, and the plasma frequency scales inversely to $d$, making the dielectric permittivity in the rods large and negative. The expressions for the series coefficients (a.k.a., dynamic correctors) and the radius of convergence in $\eta$ are explicitly related to the solutions of higher-order cell problems and the geometry of the rods. Within the radius of convergence, we are able to compute the dispersion relation and the fields and define dynamic effective properties in a mathematically rigorous manner. Explicit error estimates show that a good approximation to the true dispersion relation is obtained using only a few terms of the expansion. 
The convergence proof requires the use of properties of the Catalan numbers to show that the series coefficients are exponentially bounded in the $H^1$ Sobolev norm.
\end{abstract}
%

\section{Introduction}\label{sec:introduction}

Sub-wavelength plasmonic crystals are a class of {\em meta-material} that possesses a microstructure consisting of a periodic array of plasmonic inclusions embedded within a dielectric host.  The term ``sub-wavelength'' refers to the regime in which the period of the crystal is smaller than the wavelength of the electromagnetic radiation traveling inside  the crystal. Many recent investigations into the behavior of meta-materials  focus on phenomena associated with the {\em quasi-static limit} in which the ratio of the period cell size to wavelength tends to zero.  
Sub-wavelength micro-structured composites are known to exhibit effective electromagnetic properties that are not available in naturally-occurring materials.
Investigations over the past decade have explored a variety of meta-materials, including arrays of micro-resonators, wires, high-contrast dielectrics, and plasmonic components.  The first two, especially in combination, have been shown to give rise to unconventional bulk electromagnetic response at microwave frequencies (Smith \textit{et al.} 2000; Pendry \textit{et al.} 1999; Pendry \textit{et al.} 1998) and, more recently, at optical frequencies Povinelli (2009), including negative effective dielectric permittivity and/or negative effective magnetic permittivity.  An essential ingredient in creating this response are local resonances contained within each period due to extreme properties, such as high conductivity and capacitance in split-ring resonators Pendry \textit{et al.} (1999).

In the case of plasmonic crystals, the dielectric permittivity $\epsp$ of the inclusions is frequency dependent and negative for frequencies below the plasma frequency~$\omegap$,
\begin{equation}\label{epsp}
  \epsp(\omega) = 1 - \frac{\omegap^2}{\omega^2}.
\end{equation}
Shvets \& Urzhumov (2004, 2005) have investigated plasmonic crystals in which  $\omegap$ is inversely proportional to the period of the crystal and for which both inclusion and host materials have unit magnetic permeability. They have proposed that simultaneous negative values for both an effective $\eps$ and $\mu$ arise at sub-wavelength frequencies that are quite far from the quasi-static limit, that is, 
\begin{equation}
\eta=kd=\frac{2\pi d}{\lambda}
\end{equation}
is not very small, where $d$ is the period of the crystal, $k$ is the norm of the Bloch wavevector and $\lambda$ is the wavelength. In this work, we present rigorous analysis of this type of plasmonic crystal by establishing the existence of convergent power series in $\eta$ for the electromagnetic fields and the first branch of the associated dispersion relation. The effective permittivity and permeability defined according to Pendry \textit{et al.} (1999) are shown to be positive for all $\eta$ within the radius of convergence $R$, and, in this regime, {\em the extreme property of the plasma produces no resonance in the effective permittivity or permeability.} This regime is well distanced from the resonant regime investigated in Shvets \& Urzhumov (2004, 2005). 

The analysis shows that the radii of convergence of the power series is at least $R_m$, which is not too small, as shown in Table \ref{radiustable}, which contains values of $R_m$ for circular inclusions of various radii $rd$. 
\begin{table}
\centering 
\begin{tabular}{|c|c|c|c|c|c|}
\hline
         $r$   &$0.1$  &$0.2$   &$0.3$     &$0.4$    &$0.45$    \\ \hline
         $R_m$   &$1/60$ &$1/68$  &$1/88$    &$1/96$   &$1/340$   \\ \hline
\end{tabular}\label{radiustable}
\caption{Lower bounds on the radii of convergence $R$ for circular inclusions of radii $rd$.}
\end{table}
The number $R_m$ can be put in physical perspective by fixing the cell size and introducing the parameters $\lambda_m$ and $k_M$ such that the power series describes wave propagation for wavelengths above  $\lambda_m$  and wave numbers below $k_M$. Table 2 contains values of $\lambda_m$ and $k_M$ when $d=10^{-7}m$. The wavelengths lie in the infrared range and the plasma frequency is $\omega_p = 10^{15}\text{sec}^{-1}$.
\begin{table}
\centering 
\begin{tabular}{|c|c|c|c|c|c|}
\hline
         $r$        &$0.1$     &$0.2$     &$0.3$     &$0.4$     &$0.45$    \\ \hline
$\lambda_m$&$38\,\mu m$&$43\,\mu m$&$56\,\mu m$&$73\,\mu m$&$214\,\mu m$ \\ \hline
  $k_M$    &$1.6\cdot10^5\,m^{-1}$ &$1.4\cdot10^5\,m^{-1}$&$1.1\cdot10^6\,m^{-1}$&$1.0\cdot10^5\,m^{-1}$&$3.0\cdot10^4\,m^{-1}$ \\ \hline
\end{tabular}\label{lambdakay}
\caption{Values of $\lambda_m$ and $k_M$ for circular inclusions of radii $rd$ when $d=10^{-7}m$.}
\end{table}

We focus on harmonic H-polarized electromagnetic waves in a lossless composite medium
consisting of a periodic array of plasmonic rods embedded
in a non-magnetic frequency-independent dielectric host material. Each period can contain
multiple parallel rods with different cross-sectional shapes, however the rod-host configurations are restricted to those with rectangular symmetry, i.e., configurations invariant under a $180^o$ rotation about the center of the unit cell. The regime of interest for this
investigation is that in which
\begin{itemize}
\item[1.] the plasma frequency $\omegap$ is high
\item[2.] the ratio of the cell width to the wavelength is small ($\eta\ll 1$).
\end{itemize}
From the formula $\epsp = 1-\omegap^2/\omega^2$, it is seen that a high plasma frequency $\omegap$ gives rise to a {\em large and negative dielectric permittivity} $\epsp$ in the plasmonic inclusions. Following Shvets \& Urzhumov (2004), the plasma frequency is related to the cell size by
$$
\omegap=\frac{c}{d}.
$$
This results in the relation
$$
\epsp=1-\frac{c^2}{\omega^2d^2},
$$
where $c$ is the speed of light in vacuum. The governing family of differential equations for the magnetic field is the Helmholtz equation with a rapidly oscillating coefficient
\begin{equation}\label{Aeta}
  -\nabla\cdot\left( A(x/d) \nabla u \right) = \frac{\omega^2}{c^2} u\,,
\end{equation}
in which $A$ is the matrix defined on the unit period of the crystal by
\[A(y)= \left\{ 
\begin{array}{ll}
\epsp^{-1}I&\mbox{in the plasmonic phase, } \\
\epspbar^{-1}I&\mbox{in the host phase, } 
\end{array}
\right. \] 
$I$ is the identity matrix and
$$
\epspbar=1\,\,\,\,\,\,\text{and}\,\,\,\,\,\,\epsp=1-\frac{c^2}{\omega^2d^2}.
$$
The coefficient $A$ is not coercive in the regime $\omegap\geq\omega$, since $\epsp$ is negative in this regime, and
%
%
%
%
%
\begin{figure}
\centerline{\scalebox{0.25}{\includegraphics{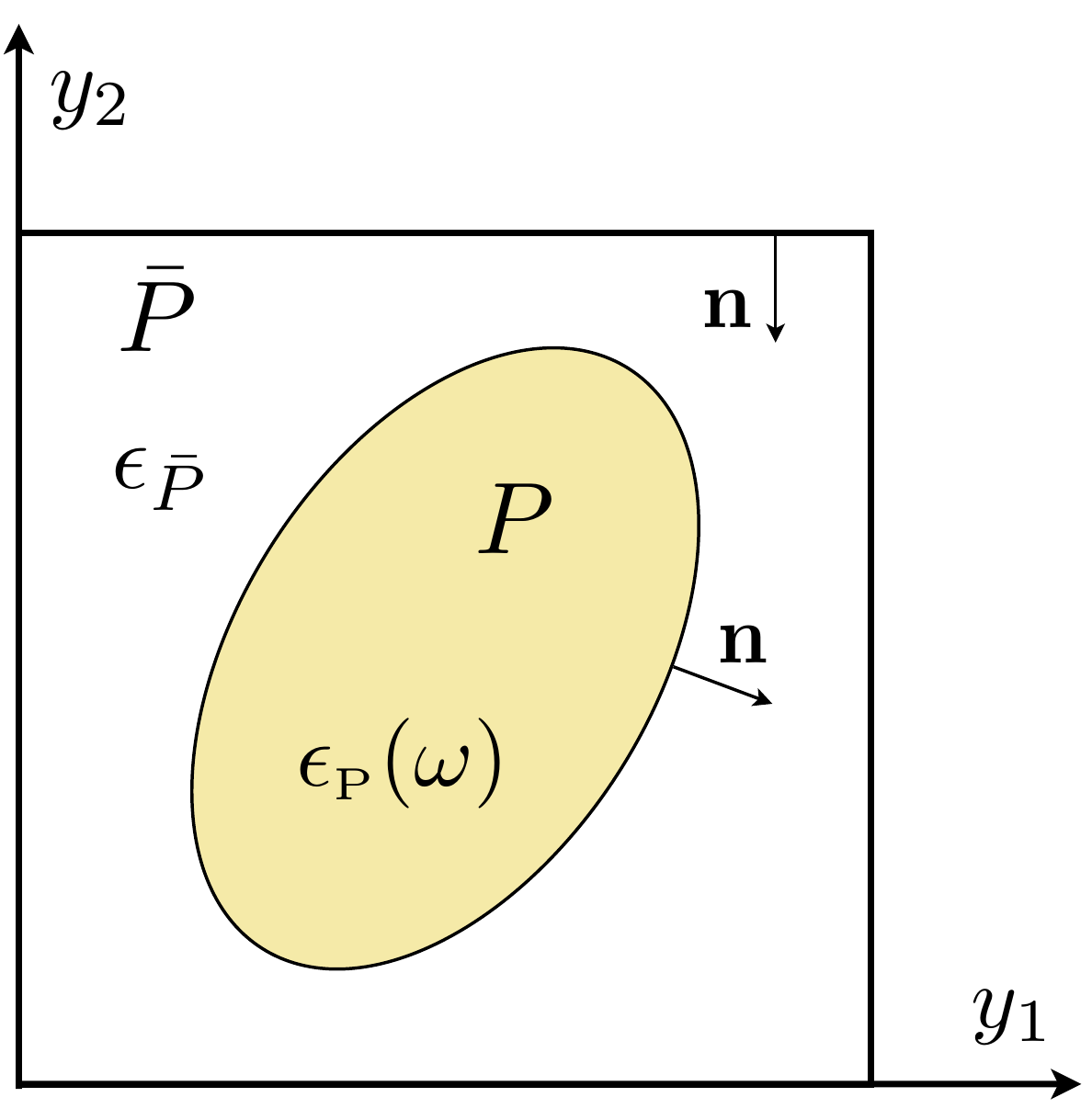}}}
\caption{Unit cell with plasmonic inclusion.}
\label{fig:inclusion}
\end{figure}
%
%
%
%
%
it is precisely the appearence of negative $\epsp$ that allows us to obtain a {\em convergent} power series expansion of the electromagnetic field and the frequency~$\omega^2$ for a fixed Bloch wavevector~$\mathbf{k}=k\kk$, with $|\kk|=1$.  

In Theorem \eqref{thm:solution}, we obtain the following series expansion for the frequency $\omega^2$
\begin{equation}\label{dispersionrelation}
  \omega^2=\omega_\eta^2=c^2k^2\sum_{m=0}^{\infty}\xi^2_{2m}\eta^{2m}
\end{equation}
in which $\xi^2_{2m}$ is a tensor of degree $2m+2$ in $\kk$. This gives rise to a convergent power series for an effective index of refraction $n^2_{\text{eff}}$ defined through
\begin{equation}\label{index}
n^2_{\text{eff}}=\frac{c^2k^2}{\omega^2_\eta}.
\end{equation}
The effective property $n^2_{\text{eff}}$ is well defined for all $\eta$ in the radius of convergence and is not phenomenological in origin but instead follows from first principles using the power series expansion. Interpreting the first term of this series as the quasi-static index of refraction $n_\text{qs}^2$, the remaining terms then provide the dynamic correctors of all orders. In section 6, we define the effective permeability $\mueff$ and prove that $n^2_{\text{eff}}$ and $\mueff$ are both positive for $\eta$ in some interval $(0,\eta_0]$ and that a mild effective magnetic response emerges for the homogenized composite, even though the component materials are non-magnetic ($\mu_{\bar P}=\mu_P=1$). Having defined $n^2_{\text{eff}}$ and $\mueff$, the effective electrical permittivity $\epseff$ can be defined through the equation
$$n^2_{\text{eff}}=\mueff\epseff,$$
so that $\epseff$ is positive whenever both $n^2_{\text{eff}}$ and $\mueff$ are positive. Thus, one has a solid basis on which to assert that plasmonic crystals function as materials of positive index of refraction in which both the effective permittivity and permeability are positive. The method developed here can be applied to other types of frequency-dependent dielectric media such as polaratonic crystals. From a physical perspective, this work provides the first explicit description of Bloch wave solutions associated with the first propagation band inside nanoscale plasmonic crystals. In the context of frequency independent dielectric  inclusions, the first two terms of $n^2_{\text{eff}}$ are identified via Rayleigh sums in McPhedran \emph{et al.} (2006). 

To emphasize the difference between effective properties defined for meta-material structures where the crystal period $d$ is fixed and effective properties defined in the quasistatic limit, i.e., $k$ fixed and $d\rightarrow 0$, we refer to the latter as quasistatic effective properties and denote these with the subscript $\text{qs}$.
The situation considered in this paper contrasts with the case in which $\eps\approx d^{-2}$ in the inclusion and is large and positive, investigated by Bouchitt\'e \& Felbacq (2004). In that case for $\eta\rightarrow 0$, $\mu_\text{qs}(\omega)$ has poles at Dirichlet eigenvalues of the inclusion and therefore is negative in certain frequency intervals (see also Bouchitt\'e \& Felbacq (2004), (2005), (2005a)). In fact, what allows us to prove convergence of the power series in the plasmonic case is precisely the absence, due to negative $\epsp$, of these internal Dirichlet resonances.

%
%
%
%
%

In the regime where $\epsp$ is negative and large, the perturbation methods used for describing Bloch waves in heterogeneous media developed in Odeh \& Keller (1964), Conca (2006), Bensoussan \textit{et al.} (1978) cannot be applied.  Our analysis instead makes use of the fact that $\epsp$ is negative and large for sub-wavelength crystals and develops high-contrast power series solutions for the nonlinear eigenvalue problem that describes the propagation of Bloch waves in plasmonic crystals.
The convergence analysis takes advantage of the iterative structure appearing in the series expansion and is inspired by a technique 
of Bruno (1991) developed for series solutions to quasi-static field problems.
We prove that the series converges to a solution of the harmonic Maxwell system for ratios of cell size to wavelength that are not too small.  Indeed for typical values of the plasma frequency $\omegap$ the analysis delivers convergent series solutions for nano scale plasmonic rods at infrared wavelengths.

In section \ref{sec:errorbounds} we compute the first two terms of the dispersion relation for circular inclusions Shvets \& Urzhumov (2004, 2005) and provide explicit bounds on the relative error comitted upon replacing the full series with its first two terms. The error is seen to be less than $3\%$ for values of $\eta$ up to $20\%$ of the convergence radius, so that the two-term approximation provides a numerically fast and accurate approximation to the dispersion relation.

The high contrast in $\eps$ gives rise to effective constants $\epseff$ and $\mueff$.  In the bulk relation
\begin{equation}\label{Beff}
  B_\text{eff} = \mueff H_\text{eff},
\end{equation}
where $B_\text{eff}$ is the average over the the period cell (a flux), whereas $H_\text{eff}$ is the average of $H_3$ over line segments in the matrix parallel to the rods.  Taking the ratio of $B_\text{eff}/H_\text{eff}$ delivers an effective magnetic permeability and one recovers magnetic activity from meta-materials made from non-magnetic materials.  This phenomenon was understood by Pendry {\it et al.} (1999), and has been made rigorous in the quasistatic limit through two-scale analysis in several cases.  These include the two-dimensional arrays of inclusions in which $\eps$ scales as $d^{-2}$ (Bouchitte \& Felbacq (2004, 2005); Felbacq \& Bouchitte (2005, 2005a); Pendry \textit{et al.} (1999)) two dimensional arrays of ring resonators whose surface conductivity scales as $d^{-1}$ Kohn \& Shipman (2008), as well as three-dimensional arrays of split-ring wire resonators in which the conductivity scales as $d^{-2}$ Bouchitte \& Schweizer (2008). This ``non-standard" homogenization has been understood for some decades in problems of porous media and imperfect interface (Cioranescu \& Paulin (1979); Auriault \& Ene (1994); Lipton (1998); Donnato \& Monsouro (2004); Zhikov (1995)) and recently has given rise to interesting effects in composites of both high contrast and high anisotropy (Cherednichenko \textit {et al.} (2006); Smyshlyaev (2009)).

The two-scale analysis in these cases relies on the coercivity of the underlying partial differential equations.  The problem of plasmonic inclusions, however, is not coercive because $\eps$ is negative in the plasma---but it is precisely this negative index that underlies the convergence of the power series.  As we shall see, the uniqueness of the solution of the Dirichlet problem for $\triangle u - u = 0$ in the plasmonic inclusion gives exponential bounds on the coefficients of the series, which allows us to prove that it converges to a solution of the differential equation \eqref{Aeta}.  This result presented here shows that by considering a finite number of terms in the series, one has an approximation of the true solution, to any desired algebraic order of convergence.  
In this context we point out the recent work of (Smyshlyaev \& Cherednichenko (2000); Kamotski \textit{et al.} (2007); Panasenko (2009)) that shows that the power series expressed by the formal two scale expansion of Bakhvalov \&  Panasenko (1989) is an asymptotic series in certain cases under the hypothesis that the coefficient $A$ is coercive.

\section{Mathematical Formulation and Background}\label{sec:formulation}

We introduce the nonlinear eigenvalue problem describing the propagation of Bloch waves inside a plasmonic crystal and provide the context for the power series approach to its solution.

For points $\xx=(x_1,x_2)$ in the $x_1x_2$-plane, the $d$-periodic dielectric coefficient of the crystal is denoted by $\eps(\omega,\xx)$, where
\[\eps(\omega,\xx)= \left\{ 
\begin{array}{lc}
\epsp(\omega)&\mbox{for } \xx\in P,\\
\epspbar&\mbox{for } \xx\in\overline{P}.
\end{array}
\right. \] 
Both materials are assumed to have unit magnetic permeability, $\mup=\mupbar=1$.

We assume a Bloch-wave form of the field, where $\hat{\kappa}=(\kappa_1,\kappa_2)$ is the unit vector along the direction of the traveling wave and $k=2\pi/\lambda$ is the wave number for a wave of length~$\lambda$. 
The magnetic and electric fields are denoted by ${\bf H}=(H_1,H_2,H_3)$ and ${\bf E}=(E_1,E_2,E_3)$ respectively.
For $H$-polarized time-harmonic waves, the non-vanishing field components are 
\begin{eqnarray}
H_3=H_3(\xx)e^{i(k\kk\cdot\xx-\omega t)},\,\,E_1=E_1(\xx)e^{i(k\kk\cdot\xx-\omega t)},\,\,
E_2=E_2(\xx)e^{i(k\kk\cdot\xx-\omega t)} \label{em3}
\end{eqnarray}
in which the fields $H_3(\xx)$, $E_1(\xx)$, and $E_2(\xx)$ are continuous and $d$-periodic in both $x_1$ and~$x_2$.  The Maxwell equations take the form of the Helmholtz equation \eqref{Aeta}, in which substitution of $u=H_3({\xx})e^{i(k\kk\cdot\xx-\omega t)}$ gives
\begin{eqnarray}
-(\nabla+ik\hat{\kappa})\epsp^{-1}(\nabla+ik\hat{\kappa})H_3&=&\frac{\omega^2}{c^2}H_3\quad\hbox{in the rods,}\label{plasmonic}\\
-(\nabla+ik\hat{\kappa})\epspbar^{-1}(\nabla+ik\hat{\kappa})H_3&=&\frac{\omega^2}{c^2}H_3\quad\hbox{in the host material,}\label{host}
\end{eqnarray}
where $H_3$ satisfies the transmission conditions on the interface between the rods and host material given by
\begin{eqnarray}
\nn\cdot(\epsp^{-1}(\nabla+ik\kk)H_3)_{|_{\scriptscriptstyle{p}}}=\nn\cdot(\epspbar^{-1}(\nabla+ik\kk)H_3)_{|_{\pbar}}.
\label{jump}
\end{eqnarray}
Here, the subscripts indicate the side of the interface where the quantities are evaluated and ${\bf n}$ is the unit normal vector to the interface pointing into the host material.  We denote the unit vector pointing along the $x_3$ direction by ${\bf e}_3$, and the electric field component of the wave is given by
${\bf E}=-\frac{ic}{\omega\epsilon}{\bf e}_3\times \nabla H_3$.

For each value of the wave-vector $\bf k$, equations (\ref{plasmonic}, \ref{host}, \ref{jump}) provide a nonlinear eigenvalue problem for the pair $H_3$ and $\omega^2$. One of the main results of this work is to show that this problem is well posed by explicitly constructing solutions using power series expansions. In order to develop the appropriate expansions, we rewrite (\ref{plasmonic}, \ref{host}, \ref{jump})  in terms of $\eta$
and a dimensionless variable $\yy$ in $\RR^2$ that normalizes a period cell to the unit square $\QQ=[0,1]^2$,
$\xx = \yy d = \yy\eta/k$. We define the $Q$-periodic function
%

$$
  h(\yy) = H_3(\yy d)
$$
and for convenience of notation, we redefine
\begin{equation}
\eps(\omega,\yy)= \left\{ 
\begin{array}{lc}
\epsp(\omega)&\mbox{for } \yy\in P,\\
\epspbar&\mbox{for } \yy\in\overline{P},
\end{array}
\right.
\label{normalizedfreq}
\end{equation}
to arrive at the eigenvalue problem that requires the pair $h(\yy)$ and $\omega^2$ to be a solution of the master system
\begin{equation}\label{master}
  \renewcommand{\arraystretch}{1.5}
\left\{
  \begin{array}{l}
    -(\grady + i\eta\kk) \eps^{-1}(\omega,\yy) (\grady + i\eta\kk)h(\yy) = \eta^2\frac{\omega^2}{c^2k^2} h(\yy)\quad\hbox{for } \yy\in P\cup\overline{P},\\
    \nn\cdot \epsp^{-1}(\omega)(\grady + i\eta\kk) h(\yy) \at{\pp} = \nn\cdot\epspbar^{-1}(\grady + i\eta\kk) h(\yy) \at{\bar P}\quad\hbox{for } \yy\in\partial P.
  \end{array}
\right.
\end{equation}
We prove in Theorem \ref{thm:solution} that this eigenvalue problem can be solved by constructing explicit convergent power series solutions. 

The development of the remainder of the paper is as follows.  In section \ref{sec:powerseries} the power series expansion is introduced  and the associated boundary-value problems necessary for determining each term in the series are obtained. The boundary value problems are given by a strongly coupled infinite system of linear partial differential equations. The existence and uniqueness of the solution to this infinite system is proved under fairly general hypotheses in section \ref{sec:solutions}.  Because the system is coupled through convolution products, the convergence analysis is delicate.  The convolutions are handled through estimates involving sequences of Catalan numbers whose convolution products determine the next element of the sequence.  The Catalan numbers and their relevant properties are discussed and used to derive bounds on the Sobolev norm of each term of the series expansion in section \ref{sec:convergence}.  These bounds are then used to establish the radius of convergence for the power series representations of the field and frequency (the dispersion relation), which solve the nonlinear eigenvalue problem \eqref{master}.  Section \ref{sec:errorbounds} deals with the computation of error bounds for finite-term approximations of the magnetic field and the dispersion relation.

\section{\bf Power Series Expansions}
\label{sec:powerseries}

We take $\eta$ to be the expansion parameter for the field $h(\yy)$ and the frequency $\omega^2$,
\begin{eqnarray}
  h_\eta        &=& h_0 + \eta h_1 + \eta^2 h_2 + \dots,\notag\\
  \omega_\eta^2 &=& \omega^2_0 + \eta\omega^2_1 + \eta^2\omega_2^2 + \dots,\label{expansions}
\end{eqnarray}
in which the functions $h_m$ are periodic with period cell $\QQ$.
%

Inserting \eqref{expansions} into \eqref{master} and identifying  
coefficients of like powers of $\eta$ on the right- and left-hand  
sides yields the equations
\begin{equation*}
  \renewcommand{\arraystretch}{1.5}
\left\{
  \begin{array}{ll}
    \Delta h_m + 2i\kdotgrad h_{m-1} - h_{m-2} = -\frac{\omega^2_\ell}
{c^2k^2} h_{m-2-\ell} & \text{in $\bar P$,} \\
    \Delta h_m + 2i\kdotgrad h_{m-1} - h_{m-2} = h_m - \frac{\omega^2_
\ell }{c^2k^2}h_{m-2-\ell} & \text{in $P$,} \\
    \left(\frac{\omega^2_\ell }{c^2k^2}\nabla h_{m-2-\ell}- \nabla  
h_m - ih_{m-1}\kk\right)\at{\bar P}\cdot\nn
              = \frac{\omega^2_\ell }{c^2k^2}\nabla h_{m-2-\ell}
\at{\pp}\cdot\nn & \text{on $\dP$,}
  \end{array}
\right.
\end{equation*}
for $m=0,1,2,\dots$, in which $h_m\equiv0$ and $\omega^2_m=0$ for $m<0$ and the terms involving the subscript $\ell$ are convolutions written according to the following summation conventions,
\begin{eqnarray*}
  a_\ell b_{n-\ell}=\sum_{\ell=0}^n a_\ell b_{n-\ell}\,, &&a_{\ell}b_{n-\ell}^{(\ell<\ell_2)}=\sum_{\ell=0}^{\ell_2-1} a_\ell b_{n-\ell},\\
   a_\ell b_{n-\ell}^{(\ell_1<\ell<\ell_2)}=\sum_{\ell=\ell_1+1}^{\ell_2-1} a_\ell b_{n-\ell}\,,&&a_\ell b_{n-\ell}^{(\ell\,\text{even})}=\sum_{\ell=0}^{[n/2]} a_{2\ell} b_{n-2\ell},
\end{eqnarray*}
where $[n/2]$ denotes the largest integer less than or equal to $n/2$. 
The boundary value problem satisfied by $h_0$ in $\bar P$ is
\begin{equation*}
  \renewcommand{\arraystretch}{1.5}
\left\{
  \begin{array}{ll}
    \Delta h_0 = 0 & \text{in $\bar P$,} \\
    \nabla h_0\at{\bar P}\cdot\nn=0 & \text{on $\dP$,}
  \end{array}
\right.
\end{equation*}
so that this function is necessarily a constant in $\bar P$. We denote this constant value by $\bar h_0$. It will be convenient to use the dimensionless parameters $\psi_m$ and $\xi^2_m$ defined through
\begin{equation}
h_m=i^m\bar h_0\psi_m\,\,\,\,\,\,\,\text{and}\,\,\,\,\,\,\omega^2_m=c^2k^2\xi^2_m.\label{dimensionless}
\end{equation}
%
%
%
In terms of $\psi_m$ and $\xi^2_m$, the above equations for the functions $h_m$ become
\begin{equation}\label{coefficientequations}
  \renewcommand{\arraystretch}{1.5}
\left\{
  \begin{array}{ll}
    \Delta \psi_m + 2\kdotgrad \psi_{m-1} + \psi_{m-2} = (-i)^\ell\xi^2_\ell \psi_{m-2-\ell} & \text{in $\bar P$,} \\
    \Delta \psi_m + 2\kdotgrad \psi_{m-1} + \psi_{m-2} = \psi_m + (-i)^\ell\xi^2_\ell \psi_{m-2-\ell} & \text{in $P$,} \\
    (\nabla((-i)^\ell\xi^2_\ell \psi_{m-2-\ell}) + \nabla \psi_m + \psi_{m-1}\kk)\at{\bar P}\cdot\nn
              = \nabla((-i)^\ell\xi^2_\ell \psi_{m-2-\ell})\at{\pp}\cdot\nn & \text{on $\dP$,}
  \end{array}
\right.
\end{equation}
for $m=0,1,2,\dots$, in which $\psi_m\equiv0$ and $\xi^2_m=0$ for $m<0$.
We thus have an infinite system which the sequences $\{\psi_m\}_{m=0}^\infty$ and $\{\xi_m^2\}_{m=0}^\infty$ must satisfy. The system is written in terms of  a Poisson equation in $\bar P$ with Neumann boundary data and a Helmholtz equation in $P$ with Dirichlet boundary data, namely
\begin{equation}\label{poisson}
  \renewcommand{\arraystretch}{1.5}
\left\{
  \begin{array}{ll}
    \Delta \psi_m =G_m & \text{in $\bar P$,} \\
     \nabla \psi_m\at{\bar P}\cdot\nn=F_m & \text{on $\dP$,}
  \end{array}
\right.
\end{equation}
and
\begin{equation}\label{helmholtz}
  \renewcommand{\arraystretch}{1.5}
\left\{
  \begin{array}{ll}
    \Delta \psi_m =\psi_m+G_m & \text{in $P$,} \\
     \psi_m\at{\pp}=\psi_m\at{\bar P} & \text{on $\dP$,}
  \end{array}
\right.
\end{equation}
in which 
\begin{equation}\label{bvpdata}
  \renewcommand{\arraystretch}{1.5}
\left\{
  \begin{array}{ll}
G_m=(-i)^\ell\xi^2_\ell \psi_{m-2-\ell}-2i\kdotgrad \psi_{m-1}-\psi_{m-2}\\
F_m=\nabla((-i)^\ell\xi^2_\ell \psi_{m-2-\ell})\at{\pp}\cdot\nn-(\nabla((-i)^\ell\xi^2_\ell \psi_{m-2-\ell})-\psi_{m-1}\kk)\at{\bar P}\cdot\nn,\\
\end{array}
\right.
\end{equation}
where $\psi_m\at{\bar P}$ and $\psi_m\at{\pp}$ indicate the trace of $\psi_m$ on the $\bar P$ and $P$ sides of the interface $\partial P$ separating the two materials. The Neumann boundary value problem \eqref{poisson} is subject to the standard \textit{solvability condition} given by
\begin{equation}\label{solvability}
\integraloverpbar{G_m}+\langle F_m\rangle_{\partial P}=0.
\end{equation}
Here the area integral over the domains $\bar P$ and $P$ are denoted by $\langle\cdot\rangle_{\bar P}$ and $\langle 
\cdot \rangle_{P}$, while the line integral over the interface $\partial P$ separating the two materials is denoted by $\langle \cdot  \rangle_{\partial P}$. 

The iterative algorithm for solving the system is as follows. First note from the definition of $\psi_0$ it follows that  that $\psi_0=1$ for $\yy$ in $\bar P$. The function $\psi_0$ is determined inside $P$ by solving (\ref{helmholtz}) with Dirichlet boundary data $ \psi_0\at{\pp}=\psi_0\at{\bar P}=1$ on $\dP$. Then $\psi_1$ on $\bar P$ is the solution of (\ref{poisson}) with Neumann data  $\nabla \psi_1\at{\bar P}\cdot\nn=-\psi_0\at{\bar P}{\hat{\kappa}}\cdot\nn$ on $\dP$. The process then continues with the boundary values on $\partial P$ of $\psi_m$ in $\bar P$ providing the Dirichlet data for $\psi_m$ in P which, in turn, provides the Neumann data for $\psi_{m+1}$ in $\bar P$, up to an additive constant. The term $\xi^2_{m-2}$ is determined by the consistency condition (\ref{solvability}) and an inductive argument can be used to show that it is a monomial of degree $m$ in $\kk$. The equations satisfied by $\psi_0,\ldots,\psi_4$ inside $\bar P$, $P$, and $\dP$ are listed in Table \ref{pdes} below.
\begin{table}\label{pdes}
$
\scriptsize{
\left.
  \begin{array}{|c|c|c|} \hline 
    \bar P & P & \dP \\ \hline\hline
    \psi_0 = 1 & \Delta\psi_0 = \psi_0 & \nabla\psi_0\at{\pbar}\cdot\nn = 0 \\ \hline
    \Delta\psi_1 + 2\kdotgrad\psi_0 = 0 & \Delta\psi_1 + 2\kdotgrad\psi_0 = \psi_1 & 
       (\nabla\psi_1 + \psi_0\kk)\at{\pbar}\cdot \nn = 0 \\ \hline
    \Delta\psi_2 + 2\kdotgrad\psi_1 + \psi_0 = \xi^2_0\psi_0 &
        \Delta\psi_2 + 2\kdotgrad\psi_1 + \psi_0 = \psi_2 + \xi^2_0\psi_0 & 
        \begin{array}{c}(\nabla(\xi^2_0\psi_0)+\nabla\psi_2+\psi_1\kk)\at{\pbar}\cdot\nn \\= \nabla(\xi^2_0\psi_0)\at{\pp}\cdot\nn\end{array} \\ \hline
    \Delta\psi_3 + 2\kdotgrad\psi_2 + \psi_1 = \xi^2_0\psi_1 & 
        \Delta\psi_3 + 2\kdotgrad\psi_2 + \psi_1 = \psi_3 + \xi^2_0\psi_1 & 
        \begin{array}{c}(\nabla(\xi^2_0\psi_1)+\nabla\psi_3+\psi_2\kk)\at{\pbar}\cdot\nn \\= \nabla(\xi^2_0\psi_1)\at{\pp}\cdot\nn\end{array} \\ \hline
    \Delta\psi_4 + 2\kdotgrad\psi_3 + \psi_2 = \xi^2_0\psi_2 - \xi^2_2\psi_0& 
        \begin{array}{c}\Delta\psi_4 + 2\kdotgrad\psi_3 + \psi_2 \\= (\xi^2_0\psi_2 - \xi^2_2\psi_0)+ \psi_4\end{array} & 
        \begin{array}{c}(\nabla(\xi^2_0\psi_2-\xi^2_2\psi_0)+\nabla\psi_4+\psi_3\kk)\at{\pbar}\cdot\nn \\= \nabla(\xi^2_0\psi_2-\xi^2_2\psi_0)\at{\pp}\cdot\nn\end{array} \\ \hline
  \end{array}
\right.}
$
\caption{Table of PDEs for the functions $\psi_m$ obtained from the expansion in $\eta$.}
\end{table}

Note in the table that $\xi^2_{odd}=0$ (meaning $\xi^2_{\ell}=0$ for $\ell=1,3,5,\dots$). In the next section we identify a large class of shapes for the plasmonic rod cross sections for which the sequences $\{\psi_m(\yy)\}_{m=0}^\infty$ and $\{\xi_m^2\}_{m=0}^\infty$  satisfy the infinite system (\ref{poisson}, \ref{helmholtz}, \ref{bvpdata}, \ref{solvability}) and $\langle \psi_m\rangle_{\bar P}=0$, $m=1,2,\ldots$. The mean zero property of $\psi_m$ on $\bar P$ provides a tractable scenario for proving the convergence of the resulting power series.  This topic is discussed further in section \ref{sec:solutions}. 

In what follows we will make use of the equivalent weak form of the infinite system.
To introduce the weak form we introduce the space of complex valued square integrable functions with square integrable derivatives $H^1(Q)$. For $u$ and $v$ in $H^1(Q)$ the inner product is defined by 
$
(u,v)_{H^1(Q)}=\left(\int_Q\,u\overline{v}\,d\yy+\int_{Q}\,\nabla u\cdot\nabla\overline{v}\,d\yy\right),
$
and the norm is given by $\Vert v\Vert_{H^1(Q)}=(v,v)^{1/2}_{H^1(Q)}$. The $H^1$ inner products and norms 
over $P$ and $\bar P$ are defined similarly.

The weak form of the infinite system is given in terms of the space $H_\text{per}^1(Q)$ of functions in $H^1(Q)$ that take the same boundary values on opposite faces of $Q$. The weak form of the system (\ref{poisson}, \ref{helmholtz}, \ref{bvpdata}, \ref{solvability}) is given by
\begin{eqnarray}
\integraloverp{[\nabla\sigma'_{m-2}+\kk\sigma'_{m-3}]\cdot\nabla\bar v -[\kdotgrad\sigma'_{m-3}-\sigma'_{m-2}-
\sigma''_{m-4}+\sigma'_{m-4}]\bar v}+\notag\\
+\integraloverpbar{[\nabla\sigma'_{m-2}+\kk\sigma'_{m-3}]\cdot\nabla\bar v -[\kdotgrad\sigma'_{m-3}-\sigma'_{m-2}-\sigma''_{m-4}+\sigma_{m-4}]\bar v}+\label{weakform} \\
+\integraloverpbar{[\nabla\psi_m+\kk\psi_{m-1}]\cdot\nabla\bar v -[\kdotgrad\psi_{m-1}+\psi_{m-2}]\bar v}=0,\notag
\end{eqnarray}
for all  $v\in H^1_{\tiny{per}}(\QQ)$, where $\sigma'_m=(-i)^\ell\xi^2_\ell \psi_{m-\ell}$ and $\sigma''_m=(-i)^\ell\psi_{m-\ell}\xi^2_{\ell-j}\xi^2_j$.
The equivalence between (\ref{poisson}, \ref{helmholtz}, \ref{bvpdata})  and the weak form follows from integration by parts and the solvability condition (\ref{solvability}) follows from (\ref{weakform}) on choosing the test function $v=1$ in (\ref{weakform}).

\section{Solutions of the Infinite System for Plasmonic Domains with Rectangular Symmetry}\label{sec:solutions}

The goal here is to identify solutions of the infinite system for which one can prove convergence of the associated power series with a minimum of effort. Looking ahead we note that the convergence proof is expedited when one can apply the Poincare inequality to the restriction of $\psi_m$ on $\bar P$ for $m$ greater than some fixed value. To this end we seek a solution ($\{\psi_m(\yy)\}_{m=0}^\infty$,   $\{\xi_m^2\}_{m=0}^\infty$) such that for $m\geq 1$ one has $\langle \psi_m\rangle_{\bar P}=0$  and the sequences satisfy  satisfy (\ref{poisson}, \ref{helmholtz}, \ref{bvpdata}, \ref{solvability}) or equivalently satisfy (\ref{weakform}).
We show that we can find such solutions for the class of plasmonic domains $P$ with rectangular symmetry. Here we suppose that the unit period cell is centered at the origin and the class of rectangular symmetric domains is given by the set of all shapes invariant under $180^o$ rotations about the origin. This class includes simply connected domains such as rectangles and ellipses as well as multiply connected domains. For these geometries and for each $m=1,2,3\ldots$ it is demonstrated that one can add an arbitrary constant to the restriction of the function $\psi_m$ on $\bar P$ with out affecting the solvability condition (\ref{solvability}). Under the assumption of rectangular symmetry for the inclusion $P$, we will show that there exists a pair ($\{\psi_m(\yy)\}_{m=0}^\infty$,   $\{\xi_m^2\}_{m=0}^\infty$) satisfying (\ref{weakform}) \textit{with the functions $\psi_m$ in the subspace} $\mathit{H^1_*(Q)\subseteq H_\text{per}^1(Q)}$ \textit{of real-valued functions with zero average in} $\mathit{\bar P}$. 

We now record the symmetries necessarily satisfied by any  solution $\psi_m\in H^1_*(Q)$ to (\ref{poisson}, \ref{helmholtz}, \ref{bvpdata}) for plasmonic domains with rectangular symmetry. We denote the dependence of $\psi_m$ on the unit vector $\kk$ by writing $\psi_m^\kk$ so that
\begin{itemize}
\item[(i)] $\psi^{-\kk}_{m}(y)=(-1)^m\psi^{\kk}_{m}(y),\,\forall y\in Q$
\item[(ii)] $\psi^{-\kk}_{m}(y)=\psi^{\kk}_{m}(-y),\,\forall y\in Q$.
\end{itemize}
Statement $(i)$ is true for inclusions of arbitrary shape, while  statement $(ii)$ is true only for inclusions with rectangular symmetry. Taken together, these statements imply that 
$
\psi^{\kk}_{m}(-y)=(-1)^m\psi^{\kk}_{m}(y),
$
so that $\psi^{\kk}_m$ is even or odd in $Q$ according as the index $m$ is even or odd. From its definition, $\psi_0\equiv 1$ in $\bar P$ and trivially satisfies the solvability condition (\ref{solvability}). The solvability of $\psi_m$ when $m\geq 1$ is proved by induction on $m$ using the weak form (\ref{weakform}). We have the following theorem
\begin{Theorem}\label{solvabilitytheorem}
For each $\kk$, there exists a sequence of functions $\{\psi_m\}_{m=1}^{\infty}$, $\psi_m\in H^1_*(Q)$, and a sequence of real numbers $\{\xi^2_m\}$, with $\xi^2_{odd}=0$, solving the weak form (\ref{weakform}) for each integer $m$.
\end{Theorem}
\begin{proof}
The proof is divided into the base case ($m=1$ and $m=2$) and the inductive step.

\vspace{0.1cm}

\noindent \underline{\bf{Base case}}:

\vspace{0.1cm}

The solvability for $\psi_1$ and $\psi_2$ can be established without the need to restrict to rectangular symmetric inclusions. This restriction will be necessary only in the inductive step. Setting $m=1$ and $v\equiv 1$ in \eqref{weakform}, we see that the left-hand side of \eqref{weakform} vanishes. This establishes the solvability for $\psi_1$. If we then take $\integraloverpbar{\psi_1}=0$, we have a solution $\psi_1\in H^1_*(Q)$. Setting $m=2$ and $v\equiv 1$ in \eqref{weakform}, we obtain
\begin{eqnarray*}
\integraloverp{\sigma'_0}+\integraloverpbar{\sigma'_0}-\integraloverpbar{\kk\cdot\nabla\psi_1+\psi_0}=0.
\end{eqnarray*}
Since $\langle\psi_0\rangle_{\QQ}>0$ (see Appendix) and $\integraloverp{\sigma'_0}+\integraloverpbar{\sigma'_0}=\xi^2_0\langle\psi_0\rangle_{\QQ}$, this is one equation in one unknown $\xi^2_0$. Solving for $\xi^2_0$ then gives
$
\xi^2_0=\langle\psi_0\rangle_{\QQ}^{-1}\integraloverpbar{\kk\cdot\nabla\psi_1+\psi_0}.
$
Choosing this value for $\xi^2_0$ and also taking $\integraloverpbar{\psi_2}=0$, we have a solution $\psi_2\in H^1_*(Q)$.

\smallskip

\noindent \underline{\bf{Inductive step}}:

\smallskip

Let $2n$ be an even positive integer and assume that (\ref{weakform}) has solutions $\psi_m\in H^1_*(Q)$ for $m=1,2,...,2n$, with $\xi^2_{m-2}\in \RR$ and $\xi^2_{odd}=0$. Then (\ref{weakform}) has solutions $\psi_{2n+1},\,\psi_{2n+2}\in H^1_*(Q)$ for $m=2n+1, 2n+2$ with $\xi^2_{2n-1}=0$ and $\xi^2_{2n}\in \RR$.
\vspace{0.5cm}

The solvability condition for $\psi_{2n+1}$ is obtained by setting $v=1$ and $m=2n+1$ in the weak form, namely
\begin{eqnarray*}
\integraloverp{[\kdotgrad\sigma'_{2n-2}-\sigma'_{2n-1}-
\sigma''_{2n-3}+\sigma'_{2n-3}]}+\\
+\integraloverpbar{[\kdotgrad\sigma'_{2n-2}-\sigma'_{2n-1}-\sigma''_{2n-3}+\sigma'_{2n-3}]}+\\
+\integraloverpbar{[\kdotgrad\psi_{2n}+\psi_{2n-1}]}=0.
\end{eqnarray*}
The hypothesis $\xi^2_{odd}=0$, $odd\leq 2n-2$, will imply that the convolutions $\sigma_m$, $m\leq 2n-2$, have the same even/odd property as the functions $\psi_m$. Indeed, writing out $\sigma_{2n-3}$, we have
$
\sigma'_{2n-3}=(-i)^{\ell}\xi^2_{\ell}\psi_{2n-3-\ell}^{(\ell\,\,\text{even})},
$
and since $2n-3-\ell$ is odd when $\ell$ is even, it follows that $\sigma'_{2n-3}$ is a linear combination of odd functions and is, therefore, an odd function. The same reasoning applies to all the other convolutions of index less than or equal to $2n-2$. Moreover, $\kdotgrad\sigma'_{2n-2}$ is an odd function, since $\sigma'_{2n-2}$ is even. Thus, all integrals in the consistency condition above vanish (for the integrals in $\Pbar$ we can also use the fact that all functions belong to $H^1_*(Q)$), except that
\begin{eqnarray*}
\integraloverp{\sigma'_{2n-1}}+\integraloverpbar{\sigma'_{2n-1}}&=&(-i)^{2n-1}\xi^2_{2n-1}\integraloverq{\psi_0}.
\end{eqnarray*}
Since $\integraloverq{\psi_0}>0$ (see Appendix), the solvability condition for $\psi_{2n+1}$ is simply $\xi^2_{2n-1}=0$. We thus take $\xi^2_{2n-1}=0$ to establish the existence of $\psi_{2n+1}=0$. Moreover, since $\psi_m$ and $\xi_{m-2}$ are real by the induction hypothesis, $0\leq m\leq 2n$, it follows that $\psi_{2n+1}$ is real-valued. Thus, taking $\integraloverpbar{\psi_{2n+1}}=0$, we have a solution $\psi_{2n+1}\in H^1_*(Q)$. Also, $\psi_{2n+1}$ is an odd function since its index is odd. We now proceed to the solvability of $\psi_{2n+2}$, namely
\begin{eqnarray*}
\integraloverp{[\kdotgrad\sigma'_{2n-1}-\sigma'_{2n-2}-
\sigma''_{2n-2}+\sigma'_{2n}]}+\\
+\integraloverpbar{[\kdotgrad\sigma'_{2n-1}-\sigma'_{2n-2}-\sigma''_{2n-2}+\sigma'_{2n}]}+\\
+\integraloverpbar{[\kdotgrad\psi_{2n+1}+\psi_{2n}]}=0.
\end{eqnarray*}
All terms in the above equation are real numbers, since we assumed $\psi_m$ and $\xi^2_{m-2}$ real for $0\leq m\leq 2n$, with $\xi^2_{odd}=0$, and we just took $\xi^2_{2n-1}=0$ and $\psi_{2n+1}$ is real-valued. Thus, this equation contains the only one undetermined term $(-i)^{2n}\xi^2_{2n}\integraloverq{\psi_0}$. Thus, we have one real equation with one real variable, so that taking $\xi^2_{2n}$ to be such as to solve this equation and also taking $\integraloverpbar{\psi_{2n+2}}=0$, we complete the proof of the inductive step.
\end{proof}

\section{\bf Proof of Convergence}
\label{sec:convergence}

In this section we show that the power series $\sum_{m=0}^{\infty}\bar p_m\eta^m$, $\sum_{m=0}^{\infty} p_m\eta^m$ and\linebreak $\sum_{m=0}^{\infty} \xi^2_m\eta^m$, where $\bar p_m=\Vert \psi_m\Vert_{H^1(\bar P)}$ and $p_m=\Vert \psi_m\Vert_{H^1(P)}$, converge and provide lower bounds on their radius of convergence. This will then be used to show that the pair $h_{\eta}=\sum_{m=0}^\infty h_m\eta ^m$ and $\omega^2_{\eta}=\sum_{m=0}^\infty\omega_m^2\eta ^m$ is a solution to \eqref{master}. In subsection \ref{thm:catalan}, we present the Catalan Bound, which is used to provide a lower bound on the radius of convergence of the power series. In subsection \ref{stabilityestimates}, we derive inequalities which bound $\bar p_m$, $p_m$ and $\xi^2_m$ in terms of lower index terms. In subsection \ref{catalanproperties}, we present the properties of the Catalan numbers relevant for bounding convolutions of the kind appearing in (\ref{poisson}) and (\ref{helmholtz}). In subsection \ref{catalanproof}, we use an inductive argument on the inequalities of subsection \ref{stabilityestimates} to prove the Catalan Bound. Finally, in subsection \ref{solutionproof} we prove that the pair $h_{\eta}$ and $\xi^2_{\eta}$ is a solution to the eigenvalue problem \eqref{master}.

\subsection{The Catalan Bound}\label{thm:catalan}
\noindent The following theorem is one of the central results of this paper
\begin{Theorem}(Catalan Bound)\label{Inductionbound}\\
For every integer $m$, we have that
\begin{equation}
 \bar p_m,\; p_m,\; |\xi^2_m|\leq \beta C_mJ^m\label{catalangrowth}
\end{equation}
in which $C_m$ is the $m^{th}$ Catalan number, $\beta=\text{max}\{\bar p_0,\,p_0,\,|\xi^2_0|\}$ and $J=\text{max}\{J_1,J_2\}$, where the numbers $J_1$ and $J_2$ are determined as follows: $J_1$ is the smallest value of $J$ such that \eqref{catalangrowth} holds for $m\leq 4$ and $J_2$ is the smallest value of $J$ for which the following polynomials $Q^*, R^*, S^*$ in the variable $J^{-1}$ are all less than unity
\begin{multline*}
Q^*=\Omega_{\bar P}[\,\,A\{2E(4)\beta J^{-2}1/3+E(4)\beta J^{-3}5/42+E(4)\beta J^{-4}1/21+E^2(4)\beta^2 J^{-4}5/42\}+\\
+2E(4)\beta J^{-2}1/3+2E(4)\beta J^{-3}5/42+E(4)\beta J^{-4}1/21+E^2(4)\beta^2 J^{-4}5/42+J^{-2}5/42\\
+2\Omega_{\bar P}(A\{2E(4)\beta J^{-3}5/42+2E(4)\beta J^{-4}1/21+E(4)\beta J^{-5}1/42+E^2(4)\beta^2 J^{-5}1/21\}+\\
+2E(4)\beta J^{-3}5/42+2E(4)\beta J^{-4}1/21+E(4)\beta J^{-5}1/42 \,+\\ +E^2(4)\beta^2 J^{-5}1/21+J^{-3}1/21+2J^{-2}5/42\,\,)],
\end{multline*}
\begin{eqnarray*}
R^*=AQ^*+E(4)\beta J^{-2}1/3+2J^{-1}1/3+J^{-2}5/42,
\end{eqnarray*}
\begin{eqnarray*}
S^*=4J\{\thetapbar Q^*+\thetap(E(4)\beta J^{-2}(1/3)+E^2(4)\beta^2J^{-3}(1/3)
+E(4)\beta J^{-3}(5/42))&&\notag\\
+\thetapbar(E(4)\beta J^{-2}(1/3)+E(4)\thetapbar\beta J^{-3}(5/42)+\thetap J^{-3}(1/21))\}+\notag\\
+\thetap\{(\,|\xi^2_0|R^*+|\xi^2_2|J^{-2}(1/7)+p_2J^{-2}(1/7)\,)+(0.7976\beta)\}.
\end{eqnarray*}
The constants $A$, $\Omega_{\bar P}$, $\beta$, $\thetapbar$, $\thetap$, $|\xi^2_0|$, $|\xi^2_2|$ and $p_2$ are determined by the particular choice of inclusion, while $E(4)=16C_2/C_5\leq 0.7619$. 
\end{Theorem}
All bounds obtained here are expressed in terms of the Catalan numbers, area fractions and geometric parameters that appear in the Poincare inequality and in an extension operator inequality. We start by listing these parameters and give the background for their description. It is known Ne$\check{c}$as (1967) that  any $H^1(\bar P)$ function $\phi$ can be extended into $P$ as an $H^1(Q)$ function $E(\phi)$ such that $E(\phi)=\phi$ for $\yy$ in $\bar P$ and
\begin{eqnarray}
\Vert E(\phi)\Vert_{H^1(P)}\leq A \Vert \phi\Vert_{H^1(\bar P)}
\label{extension}
\end{eqnarray}
where $A$ is a nonnegative constant and is independent of $\phi$ depending only on $P$. 
For general shapes $A$ can be calculated via numerical solution of a suitable eigenvalue problem. Constants of this type appear in Bruno (1991) for high contrast expansions of the DC fields inside frequency independent dielectric media. The second constant is the Poincare constant $D^2_{\bar P}$ given by the reciprocal of the first nonzero Neumann eigenvalue of $\bar P$ and we have that $\Omega_{\bar P}=1+D_{\bar P}^2$. The last two geometric constants appearing in the bounds are the volume fractions $\theta_{P}$ and $\theta_{\bar P}$ of the regions $P$ and $\bar P$. Using that $C_m\leq 4^m$ (see section \ref{catalanproperties}), theorem \eqref{Inductionbound} shows that $\sum\bar p_m\eta^m$, $\sum p_m\eta^m$ and $\sum\xi^2_m\eta^m$ are convergent for $\eta\leq 1/4J$, so that one may prove the following theorem
%
%
%
%
\begin{Theorem}(Solution of the Eigenvalue Problem)\label{thm:solution}
Let $R=1/4J$, where $J$ is the number prescribed by theorem \eqref{Inductionbound}. Then $\sum_{m=0}^\infty\xi^2_m\eta^m$ converges as a series of real numbers and $\sum_{m=0}^\infty \psi_m\eta^m$ converges in the $H^1(Q)$ Sobolev norm for $\eta\leq R$. Moreover, $\omega^2_\eta=c^2k^2\sum_{m=0}^\infty\xi^2_m\eta^m$ and $h_\eta=\bar h_0\sum_{m=0}^\infty i^m\psi_m\eta^m$ satisfy the eigenvalue problem given by the master system \eqref{master} (or by \eqref{Aeta}).
\label{Eigenvalue}
\end{Theorem}
%
%
%
%
%

\subsection{The $\bar p_m$, $p_m$ and $\xi^2_m$ Inequalities---Stability Estimates}\label{stabilityestimates} 
We now derive the inequalities which bound $\bar p_m$, $p_m$ and $\xi^2_m$ in terms of lower index terms. These inequalities follow from stability estimates for (\ref{poisson}, \ref{helmholtz}, \ref{bvpdata}).
\begin{Theorem}
Let $m\geq 0$ be an integer. Then 
\begin{eqnarray}
\bar p_m&\leq&\Omega_{\bar P}[A\{2q\,'_{m-2}+2q\,'_{m-3}+q\,'_{m-4}+q\,''_{m-4}\}+2\bar q\,'_{m-2}+2\bar q\,'_{m-3}+\bar    q\,'_{m-4}\notag\\
     & &+\bar q\,''_{m-4}+\bar p_{m-2} +2\,\Omega_{\bar P}(A\{2q\,'_{m-3}+2q\,'_{m-4}+q\,'_{m-5}+q\,''_{m-5}\}\notag\\
     & &+2\bar q\,'_{m-3}+2\bar q\,'_{m-4}+\bar q\,'_{m-5}+q\,''_{m-5}+\bar p_{m-3}+2\bar p_{m-2})]\notag\\
p_m &\leq& A\bar p_m+q\,'_{m-2}+2p_{m-1}+p_{m-2}\label{threeinequalities}\\
|\xi^2_{m-1}|&\leq&\integraloverq{\psi_0}^{-1}\{\sqrt{\theta_P}\,q\,'^*_{m-1}+\sqrt{\theta_{\bar P}}\bar p_{m}+\notag\\  
& & +\sqrt{\theta_P}(q\,'_{m-2}+q\,''_{m-3}+q\,'_{m-3})+\sqrt{\theta_{\bar P}}(\bar q\,'_{m-2}+\bar q\,''_{m-3}+\bar q\,'_{m-3})\},\notag
\end{eqnarray}
where the $\bar p_m$ inequality holds for $m\geq 2$ only. 
\end{Theorem}
\noindent Here we have introduced the notation
\begin{eqnarray*}
q\,'_m=|\xi^2_{\ell}|p_{m-\ell},&q\,''_m=p_{m-\ell}|\xi^2_{\ell-j}||\xi^2_j|,&q\,'^*_{m-1}=|\xi^2_{\ell}|p_{m-1-\ell}^{(\ell<m-1)}\\
\bar q\,'_m=|\xi^2_{\ell}|\bar p_{m-\ell},&\bar q\,''_m=\bar p_{m-\ell}|\xi^2_{\ell-j}||\xi^2_j|.
\end{eqnarray*}
\begin{proof}
We start by proving the $p_m$ inequality. Recalling that (\ref{helmholtz}) is satisfied by $\psi_m$ in $P$ gives
\[\left\{ 
\begin{array}{ll}
\Delta \psi_m=\psi_m +G_m, & \text{in $P$}\\
\psi_m\at{\pp}=\psi_m\at{\bar P},& \text{on $\dP$}\end{array} 
\right. \] 
where $G_m=(-i)^\ell\xi^2_\ell \psi_{m-2-\ell} -2\kdotgrad \psi_{m-1} -\psi_{m-2}$. Write the orthogonal decomposition $\psi_m=u_m+v_m$, where

\begin{equation}
\label{h1}
\left\{ 
\begin{array}{ll}
\Delta u_m=u_m, & \text{in $P$}\\
u_m=\psi_m,& \text{on $\dP$}\end{array} 
\right. 
\end{equation}
and
\[\left\{ 
\begin{array}{ll}
\Delta v_m=v_m+G_m, & \text{in $P$}\\
v_m=0,& \text{on $\dP$}.\end{array} 
\right. \]
We then have by the triangle inequality that 
\begin{equation}\label{triangleinequality}
\honepnorm{\psi_m}\leq \honepnorm{u_m}+\honepnorm{v_m}.
\end{equation} 
The term $\honepnorm{u_m}$ is bounded using (\ref{extension}) and 
\begin{eqnarray}
\honepnorm{u_m}\leq\Vert E(\psi_m)\Vert_{H^1(P)}
\label{extreme}
\end{eqnarray}
 to obtain
\begin{equation}\label{tracetheorem}
\honepnorm{u_m}\leq A \honepbarnorm{\psi_m}.
\end{equation}
Here (\ref{extreme}) follows from the fact that the solution of (\ref{h1}) minimizes the $H^1(P)$ norm
over all functions with the same trace on $\partial P$.
The term $\honepnorm{v_m}$ can be bounded using a direct integration by parts on the BVP for $v_m$ 
\begin{equation}\label{vmbound}
\honepnorm{v_m}\leq \ltwopnorm{G_m}.
\end{equation}
Now,
\begin{eqnarray}
\ltwopnorm{G_m}&=&\ltwopnorm{(-i)^\ell\xi^2_\ell \psi_{m-2-\ell} -2\kdotgrad \psi_{m-1} -\psi_{m-2}}\notag\\
      &\leq&|\xi^2_{\ell}|\ltwopnorm{\psi_{m-2-\ell}}+2|\kk|\ltwopnorm{|\nabla\psi_{m-2}|^2}+\ltwopnorm{\psi^2_{m-2}}\notag\\
      &\leq&|\xi^2_\ell| p_{m-2-\ell}+2p_{m-1}+p_{m-2},\notag
\end{eqnarray}
where $p_m=\honepnorm{\psi_m}$. Using \eqref{tracetheorem} and \eqref{vmbound} in \eqref{triangleinequality} gives
\begin{equation}
p_m\leq A\bar p_m+q\,'_{m-2}+2p_{m-1}+p_{m-2},\label{secondstability2}
\end{equation}
and the $p_m$ inequality is established. We now prove the $\bar p_m$ inequality. In the weak form (\ref{weakform}), set $v=\psi_m$ in $\bar P$ and $v=u_m$ in $P$ to obtain
\begin{eqnarray*}
\integraloverpbar{[\nabla\sigma'_{m-2}+\kk\sigma'_{m-3}]\cdot\nabla\psi_m -[\kdotgrad\sigma'_{m-3}-\sigma'_{m-2}-
\sigma''_{m-4}+\sigma'_{m-4}]\psi_m}+\\
+\integraloverp{[\nabla\sigma'_{m-2}+\kk\sigma'_{m-3}]\cdot\nabla u_m -[\kdotgrad\sigma'_{m-3}-\sigma'_{m-2}-\sigma''_{m-4}+\sigma'_{m-4}] u_m}+\\
+\integraloverpbar{[\nabla\psi_m+\kk\psi_{m-1}]\cdot\nabla\psi_m -[\kdotgrad\psi_{m-1}+\psi_{m-2}]\psi_m}=0,
\end{eqnarray*} 
We now use the Cauchy-Schwarz inequality on the product of integrals appearing in each individual term. For the convolutions, we obtain
\begin{eqnarray*}
|\vectorintegraloverp{\nabla\sigma'_{m-2}}{\nabla u_m}|&=&|\vectorintegraloverp{\nabla((-i)^{\ell}\psi_{m-2-\ell}\xi^2_{\ell})}{\nabla u_m}|\\                                                   &=&|(-i)^{\ell}\xi^2_{\ell}\vectorintegraloverp{\nabla\psi_{m-2-\ell}}{\nabla u_m}|\\
            &\leq&|\xi^2_{\ell}|\bar p_{m-2-\ell}A\bar p_m,\\
            &=&q\,'_{m-2}A\bar p_m
\end{eqnarray*}
where we used that $\honepnorm{u_m}\leq A \bar p_m$.
For the double-convolutions, we obtain
\begin{eqnarray*}
|\integraloverp{\sigma''_{m-4}u_m}|
             \leq q\,''_{m-2}A\bar p_m.
\end{eqnarray*}
Proceeding similarly with the other terms, we obtain
\begin{eqnarray}
\vectorintegraloverpbar{\nabla\psi_{m}}{\nabla\psi_{m}}&\leq &
\bar p_m(\,A\{2q\,'_{m-2}+2q\,'_{m-3}+q\,'_{m-4}+q\,''_{m-4}\}+\notag\\
& &+2\bar q\,'_{m-2}+2\bar q\,'_{m-3}+\bar q\,'_{m-4}+\bar q\,''_{m-4}+2\bar p_{m-1}+\bar p_{m-2})\label{pbargradnorm}.
\end{eqnarray}
Since the functions $\psi_m$ have zero average in $\bar P$, we have the Poincare inequality
\begin{equation}\label{poincare}
 \integraloverpbar{\psi^2_m}\leq D^2_{\bar P}\vectorintegraloverpbar{\nabla\psi_m}{\nabla\psi_m},
\end{equation}
where the constant $D_{\bar P}$ can be computed from the Rayleigh quotient characterization of the first positive eigenvalue for the free membrane problem in $\bar P$. A simple computation using \eqref{poincare} then gives
$
\bar p^2_m\leq \Omega_{\bar P}\vectorintegraloverpbar{\nabla\psi_m}{\nabla\psi_m},
$
where $\Omega_{\bar P}=D^2_{\bar P}+1$. Using this inequality (\ref{pbargradnorm}) gives:
\begin{eqnarray}
\bar p_m&\leq &\Omega_{\bar P}(\,A\{2q\,'_{m-2}+2q\,'_{m-3}+q\,'_{m-4}+q\,''_{m-4}\}+\notag\\
& &+2\bar q\,'_{m-2}+2\bar q\,'_{m-3}+\bar q\,'_{m-4}+\bar q\,''_{m-4}+\bar p_{m-2}+2\bar p_{m-1})\label{oak}.
\end{eqnarray}
It will turn out to be to our advantage to apply \eqref{oak} to the last term $2\bar p_{m-1}$ in \eqref{oak} so as to replace it with
\begin{eqnarray}
\bar p_{m-1}&\leq &\Omega_{\bar P}(\,A\{2q\,'_{m-3}+2q\,'_{m-4}+q\,'_{m-5}+q\,''_{m-5}\}+\notag\\
& &+\bar q\,'_{m-3}+2q\,'_{m-4}+\bar q\,'_{m-5}+\bar q\,''_{m-5}+\bar p_{m-3}+2\bar p_{m-2})\label{basswood}.
\end{eqnarray}
Using \eqref{basswood} in \eqref{oak} yields the $\bar p_m$ inequality in \eqref{threeinequalities}, valid for $m\geq 2$ (for $m=1$, use \eqref{oak}):

Last we establish the $\xi^2_{m-1}$ inequality. Setting $v=1$ in the weak form (\ref{weakform}) we obtain
\begin{eqnarray}
\integraloverp{\kdotgrad\sigma'_{m-3}-\sigma'_{m-2}-
\sigma''_{m-4}+\sigma'_{m-4}}+\notag\\
+\integraloverpbar{\kdotgrad\sigma'_{m-3}-\sigma'_{m-2}-\sigma''_{m-4}+\sigma'_{m-4}}+\notag\\
+\integraloverpbar{\kdotgrad\psi_{m-1}+\psi_{m-2}}=0.\label{consistency}
\end{eqnarray}
(recall that for $m$ odd, each term on the left-hand side of the above equation vanishes individually). Solving for $\xi^2_{m-2}$ we then obtain
\begin{eqnarray}5
(-i)^{m-2}\xi^2_{m-2}\integraloverq{\psi_0}=\integraloverp{\kdotgrad\sigma'_{m-3}-\sigma'^{\,\,*}_{m-2}-
\sigma''_{m-4}+\sigma'_{m-4}}+\notag\\
+\integraloverpbar{\kdotgrad\sigma'_{m-3}-\sigma'^{\,\,*}_{m-2}-\sigma''_{m-4}+\sigma'_{m-4}}+\label{consistencyagain}\\
+\integraloverpbar{\kdotgrad\psi_{m-1}+\psi_{m-2}},\notag
\end{eqnarray}
where $\sigma'^{\,\,*}_{m-2}+(-i)^{m-2}\xi^2_{m-2}=\sigma'_{m-2}$. We shall be using this equality for $m\geq 5$ only, so that $\integraloverpbar{\sigma'^{\,\,*}_{m-2}}=0$ and $\integraloverpbar{\psi_{m-2}}=0$. Moreover, using that 
$
\integraloverp{|\psi_m|}\leq \sqrt{\theta_{P}\integraloverp{|\psi_m|^2}} 
$
and 
$
\integraloverpbar{|\psi_m|}\leq \sqrt{\theta_{\bar P}\integraloverpbar{|\psi_m|^2}}, 
$
where $\theta_P$ and $\theta_{\bar P}$ denote the volume fractions of the regions $P$ and $\bar P$, we have that $\integraloverp{\psi_m}\leq \sqrt{\theta_P}\,p_m$ and $\integraloverpbar{\psi_m}\leq \sqrt{\theta_{\bar P}}\,\bar p_m$. Thus, proceeding with \eqref{consistencyagain} as we did in the previous stability estimates, we obtain 
\begin{multline*}
|\xi^2_{m-2}|\leq\integraloverq{\psi_0}^{-1}\{\sqrt{\theta_P}\,p_{m-2-\ell}|\xi^2_{\ell}|^{\ell<m-2}+\sqrt{\theta_{\bar P}}\bar p_{m-1}+\sqrt{\theta_P}(q\,'_{m-3}+q\,''_{m-4}+q\,'_{m-4})+\notag\\
+\sqrt{\theta_{\bar P}}(\bar q\,'_{m-3}+\bar q\,''_{m-4}+\bar q\,'_{m-4})\}.
\end{multline*}
Since the iteration scheme at each step involves $p_m$ and $\overline{p}_m$  and $\xi_{m-1}^2$ we adjust subscripts in the above inequality to obtain the $\xi^2_{m-1}$ inequality
\begin{eqnarray}
|\xi^2_{m-1}|&\leq&\integraloverq{\psi_0}^{-1}\{\sqrt{\theta_P}\,q\,'^*_{m-1}+\sqrt{\theta_{\bar P}}\bar p_{m}+\sqrt{\theta_P}(q\,'_{m-2}+q\,''_{m-3}+q\,'_{m-3})+\notag\\
& &+\sqrt{\theta_{\bar P}}(\bar q\,'_{m-2}+\bar q\,''_{m-3}+\bar q\,'_{m-3})\}.\label{thirdstability3}
\end{eqnarray}
\end{proof}
\subsection{The Catalan Numbers}\label{catalanproperties}
We briefly present some facts about the Catalan numbers which will be used in the sequel. The Catalan numbers $C_m$ are defined algebraically through the recursion
\begin{eqnarray}
C_{m+1}=C_{m-\ell}C_\ell,\,\,\,\,C_0=1.
\label{Cm}
\end{eqnarray}
These numbers arise in many combinatorial contexts Lando (2003) as well as in the study of fluctuations in coin tossing and random walks Feller (1968). It can be shown that $C_m=\frac{1}{m+1}\binom{2m}{m}$ and a simple computation then gives the ratio of successive Catalan numbers
\begin{equation}\label{catalanratio}
\frac{C_{m+1}}{C_m}=4-\frac{6}{m+2}\,\,\,\,\,\,\text{and}\,\,\,\,\,\,\frac{C_{m}}{C_{m+1}}=\frac{1}{4}+\frac{3}{8m+4}.
\end{equation}
The first inequality above provides the exponential bound 
\begin{eqnarray}
C_m\leq 4^m. \label{upperCmbound}
\end{eqnarray}
It will be convenient to introduce the notation $\rho_m^k=C_{m-k}/C_m$. From \eqref{catalanratio}, it is clear that $\rho_m^k$ is decreasing in both $m$ and $k$. In section \ref{catalanproof} we shall make use of Table 4, in which values of $\rho_5^k$ are listed.
\begin{table}
\centering 
\begin{tabular}{|c|c|c|c|c|c|}
\hline
         $k$      &0  &1  &2   &3   &4    \\ \hline
        $\rho_5^k$&1   &1/3  &5/42  &1/21  &1/42 \\ \hline
\end{tabular}\label{rhotable}
\caption{Values of $\rho_5^k$, where $\rho_m^k=C_{m-k}/C_m$.}
\end{table}

\subsubsection{The Even Part of the Catalan Convolution}
The fact that $\xi^2_{odd}=0$ needs to be taken into account in order to provide a suitable upper estimate on the incomplete convolution term $q\,'^*_{m-1}$ appearing in the $\xi^2_{m-1}$ inequality in \eqref{threeinequalities}. Thus we consider the convolution $C_{n-\ell}C_{\ell}$ with the odd values of the index $\ell$ omitted and denote it by $C_{n-\ell}C_{\ell}^{(\ell\,\,even)}$. We then define the even part $E(n)$ by
\begin{eqnarray}
E(n)=\frac{C_{n-\ell}C_{\ell}^{(\ell\,\,even)}}{C_{n-\ell}C_{\ell}}.
\label{oddeven}
\end{eqnarray}
The following lemma gives the estimate $E(n)\leq E(4)$, $n\geq 4$.
\begin{Lemma}\label{oddomission}
The following two statements are true for all $m\geq 0$: (i) $E(2m)$ is a decreasing sequence; and (ii) $E(2m+1)=1/2$. Thus, for all $m\geq 4$, we have that $E(m)\leq \text{max}\{E(4),1/2\}=E(4)$.
\end{Lemma}
\begin{proof}
Statement $(ii)$ is actually just an observation, as one can see by writing out the sum $C_{n-\ell}C_{\ell}$. Statement (i) can be deduced from the identity $C_{2m-2\ell}C_{2\ell}\!\!=\!\!~4^mC_m$ Koshy (2008). Indeed, dividing both sides of this identity by $C_{2m+1}$, we obtain
\begin{eqnarray*}
 E(2m)&=&\frac{1}{4}\frac{4C_m}{C_{m+1}}\frac{4C_{m+1}}{C_{m+2}}\cdots\frac{4C_{2m}}{C_{2m+1}}.
\end{eqnarray*}
From \ref{catalanratio}, each of the above fractions $4C_m/C_{m+\ell}$, $\ell=1,2,...,m+1$, is a decreasing sequence in $m$ so that their product is also decreasing in $m$. This completes the proof.
\end{proof}

\subsection{Proof of the Catalan Bound}\label{catalanproof}

\begin{proof}(Catalan Bound, theorem \eqref{Inductionbound})
Fix the values of the geometric parameters $A$, $\Omegap$, $\theta_P$ and $\theta_{\bar P}$ in \eqref{threeinequalities}. Starting with the initial estimates $\bar p_0=\theta_{\Pbar},\,p_0\leq\thetap\text{ and }\xizero\leq 1$, the inequalities \eqref{threeinequalities} can be used recursively for $m=1,\,2,\,3,\,4$ to determine a number $J_1$ such that 
\begin{equation*}
\bar p_m,\,p_m,\,|\xi^2_m|\leq \beta C_mJ_1^m,\,\,\,\,\,\,\,0\leq m\leq4.
\end{equation*}
We now proceed inductively: assume that
\begin{equation}\label{catalanhypothesis}
\bar p_n,\,p_n,\,|\xi^2_n|\leq \beta C_nJ^n,\,\,\,\,\,\,\,n\in \{0,1,2,\ldots,m-1\},
\end{equation}
where $m\geq 5$. We then get for the single convolutions
\begin{eqnarray*}
q\,'_{m-k}&=&p_{m-k-\ell}|\xi^2_{\ell}|\notag\\
          &\leq&(\beta C_{m-k-\ell}J^{m-k-\ell})(\beta C_{\ell}J^{\ell})^{(\ell\,\,even)}\notag\\
          &=&\beta^2 J^{m-k}C_{m-k-\ell}C_{\ell}^{(\ell\,\,even)}\notag\\
          &\leq&E(4)\beta^2 J^{m-k}C_{m-k-\ell}C_{\ell}\notag\\
  	  &=&E(4)\beta J^{-k}\left(\frac{C_{m+1-k}}{C_m}\right)\,\beta J^mC_m\notag\\
  	  &=&E(4)\beta J^{-k}\rho_m^{k-1}\,\beta J^mC_m.
\end{eqnarray*}
where $\rho_m^k=C_{m-k}/C_m$ and lemma \eqref{oddomission} was used to introduce the factor $E(4)$. Similarly, for double convolutions we get
\begin{eqnarray*}
q\,''_{m-k}
	   &\leq&E^2(4)\beta^2 J^{-k}\rho_m^{k-2}\,\beta J^mC_m\label{doubleconvolutions}
\end{eqnarray*}
where the factor $E(4)^2$ comes from using lemma \eqref{oddomission} twice. For the non-convolution terms we get
$p_{m-k}\leq
						J^{-k}\rho_m^k\,\beta J^mC_m.\label{nonconvolution}
$
The same bounds hold for the terms $\bar p_{m-k}$, $\bar q\,'_{m-k}$ and $\bar q\,''_{m-k}$, so that we have

\begin{eqnarray}
\bar p_{m-k},\,p_{m-k}&\leq&J^{-k}\rho_m^k\,\beta J^mC_m\notag\\
\bar q\,'_{m-k},\, q\,'_{m-k}&\leq& E(4)\beta J^{-k}\rho_m^{k-1}\,\beta J^mC_m\label{bounds}\\
\bar q\,''_{m-k},\, q\,''_{m-k}&\leq& E^2(4)\beta^2 J^{-k}\rho_m^{k-2}\,\beta J^mC_m.\notag
\end{eqnarray}

The proof now essentially consists of applying these bounds to all terms in inequalities \eqref{threeinequalities}. The factor $J^{-k}$ appearing on the right-hand side of each inequality is the workhorse of the proof: by taking $J$ sufficiently large, it will allow us to close the induction argument. The incomplete convolution term $q\,'^*_{m-1}$ presents special difficulties, since attempting a bound of the type \eqref{bounds} for this term does not produce a factor of $J^{-k}$ (actually, it produces $J^0=1$). 

Recall the $\bar p_m$ inequality from \eqref{threeinequalities}
\begin{multline*}
\bar p_m\leq\Omega_{\bar P}[A\{2q\,'_{m-2}+2q\,'_{m-3}+q\,'_{m-4}+q\,''_{m-4}\}\,+ \\
     +2\bar q\,'_{m-2}+2\bar q\,'_{m-3}+\bar    q\,'_{m-4}+\bar q\,''_{m-4}+\bar p_{m-2}\,+\\
     +2\,\Omega_{\bar P}(A\{2q\,'_{m-3}+2q\,'_{m-4}+q\,'_{m-5}+q\,''_{m-5}\}\,+\\
     +2\bar q\,'_{m-3}+2\bar q\,'_{m-4}+\bar q\,'_{m-5}+q\,''_{m-5}+\bar p_{m-3}+2\bar p_{m-2})]\notag
\end{multline*}     
Using \eqref{bounds} on this inequality gives
\begin{equation}
\bar p_m\leq Q_m\,\beta J^mC_m,\label{pbarpolynomial}
\end{equation}
where $Q_m$ is the following polynomial in $J^{-1}$
\begin{multline*}
Q_m=\Omega_{\bar P}[\,\,A\{2E(4)\beta J^{-2}\rho_m^1+E(4)\beta J^{-3}\rho_m^2+E(4)\beta J^{-4}\rho_m^3+E^2(4)\beta^2 J^{-4}\rho_m^2\}+\notag\\
+2E(4)\beta J^{-2}\rho_m^1+2E(4)\beta J^{-3}\rho_m^2+E(4)\beta J^{-4}\rho_m^3+E^2(4)\beta^2 J^{-4}\rho_m^2+J^{-2}\rho_m^2\,+\notag\\
+2\Omega_{\bar P}(A\{2E(4)\beta J^{-3}\rho_m^2+2E(4)\beta J^{-4}\rho_m^3+E(4)\beta J^{-5}\rho_m^4+E^2(4)\beta^2 J^{-5}\rho_m^3\}+\notag\\ 
+2E(4)\beta J^{-3}\rho_m^2+2E(4)\beta J^{-4}\rho_m^3+E(4)\beta J^{-5}\rho_m^4+E^2(4)\beta^2 J^{-5}\rho_m^3+J^{-3}\rho_m^3+2J^{-2}\rho_m^2\,\,)].
\end{multline*}
Since we shall be using this inequality for $m\geq 5$ only, table \ref{rhotable} can be used to bound the numbers $\rho_m^k$, so that we may write $Q_m\leq Q^*$, where
\begin{multline}
Q^*=\Omega_{\bar P}[\,\,A\{2E(4)\beta J^{-2}1/3+E(4)\beta J^{-3}5/42+E(4)\beta J^{-4}1/21+E^2(4)\beta^2 J^{-4}5/42\}+\\
+2E(4)\beta J^{-2}1/3+2E(4)\beta J^{-3}5/42+E(4)\beta J^{-4}1/21+E^2(4)\beta^2 J^{-4}5/42+J^{-2}5/42\\
+2\Omega_{\bar P}(A\{2E(4)\beta J^{-3}5/42+2E(4)\beta J^{-4}1/21+E(4)\beta J^{-5}1/42+E^2(4)\beta^2 J^{-5}1/21\}+\\
+2E(4)\beta J^{-3}5/42+2E(4)\beta J^{-4}1/21+E(4)\beta J^{-5}1/42+\\+E^2(4)\beta^2 J^{-5}1/21+J^{-3}1/21+2J^{-2}5/42\,\,)].\label{qstar}
\end{multline}
The strategy now is to determine similar polynomials $R_m$ and $S_{m-1}$ for the other two inequalities, that is
$p_m\leq R_m \,\beta J^mC_m$ and $|\xi^2_{m-1}|\leq S_{m-1} \,\beta J^{m-1}C_{m-1}$, and then take $J$ large enough that all three polynomials are less than unity, allowing us to complete the induction argument. Having obtained $Q_m$, it is straightforward to obtain $R_m$. Indeed, using \eqref{bounds} and \eqref{pbarpolynomial}, the $p_m$ inequality in \eqref{threeinequalities} yields $p_m\leq R_m\,\beta J^mC_m$, where
\begin{eqnarray*}
R_m=AQ_m+E(4)\beta J^{-2}\rho_m^1+2J^{-1}\rho_m^1+J^{-2}\rho_m^2.
\end{eqnarray*}
Thus, $R_m\leq R^*$, where
\begin{eqnarray}
R^*=AQ^*+E(4)\beta J^{-2}1/3+2J^{-1}1/3+J^{-2}5/42.\label{rstar}
\end{eqnarray}
The $\xi^2_{m-1}$ inequality requires a little more care due to the presence of the incomplete convolution term $q\,'^*_{m-1}$. For the remaining terms, we proceed as we did with the previous inequalities:
\begin{eqnarray*}
\thetapbar\bar p_m+\thetap(q\,'_{m-2}+q\,''_{m-3}+q\,'_{m-3})+\thetapbar(\bar q\,'_{m-2}+\thetapbar |\xi^2_{\ell}\xi^2_{m-3-\ell}|+\thetap |\xi^2_{m-3}|)&\leq& \\
\{\thetapbar Q_m+\thetap(0.5\beta J^{-2}\rho_m^1+0.25\beta^2J^{-3}\rho_m^1
+0.5\beta J^{-3}\rho_m^2)&&\\
+\thetapbar(0.5\beta J^{-2}\rho_m^1+0.5\thetapbar\beta J^{-3}\rho_m^2+J^{-3}\rho_m^3)\}\,\beta J^mC_m.&&
\end{eqnarray*}
since this is an upper bound on $|\xi^2_{m-1}|$, we must replace the term $\beta J^mC_m$ with $\beta J^{m-1}C_{m-1}$ as follows:
\begin{eqnarray*}
\beta J^mC_m&=   &J\left(\frac{C_m}{C_{m-1}}\right)\,\beta J^{m-1}C_{m-1}\\ 
            &\leq&4J\,\beta J^{m-1}C_{m-1}.
\end{eqnarray*}
Using this replacement and the bounds \ref{rhotable} on the numbers $\rho_m^k$, we obtain the upper bound
\begin{eqnarray}
4J\{\thetapbar Q^*+\thetap(E(4)\beta J^{-2}(1/3)+E^2(4)\beta^2J^{-3}(1/3)
+E(4)\beta J^{-3}(5/42))&&\notag\\
+\thetapbar(E(4)\beta J^{-2}(1/3)+E(4)\thetapbar\beta J^{-3}(5/42)+\thetap J^{-3}(1/21))\}\,\beta J^{m-1}C_{m-1}
\label{xione}
\end{eqnarray}
It remains to deal with $q\,'^*_{m-1}$. To do this, we first write $q\,'^*_{m-1}$ as follows
\begin{eqnarray}
q\,'^*_{m-1}=p_{m-1}|\xi^2_0|+p_{m-3}|\xi^2_2|+p_2|\xi^2_{m-3}|+p_{m-1-\ell}|\xi^2_{\ell}|^{2<\ell<m-3}.\label{sum}
\end{eqnarray}
The non-convolution terms then give
\begin{multline}
p_{m-1}|\xi^2_0|+p_{m-3}|\xi^2_2|+p_2|\xi^2_{m-3}| \\
\leq\left(|\xi^2_0|R_{m-1}+|\xi^2_2|J^{-2}\frac{C_{m-3}}{C_{m-1}}+p_2J^{-2}\frac{C_{m-3}}{C_{m-1}}\right)\,\beta J^{m-1}C_{m-1}\\
\leq (\,|\xi^2_0|R^*+|\xi^2_2|J^{-2}(1/7)+p_2J^{-2}(1/7)\,)\,\beta J^{m-1}C_{m-1},\label{xitwo}
\end{multline}
since $C_{m-3}/C_{m-1}\leq C_2/C_4=1/7$, if $m\geq 5$. The remaining term $p_{m-1-\ell}|\xi^2_{\ell}|^{2<\ell<m-3}$ is treated in a completely different manner: 
\begin{eqnarray}
p_{m-1-\ell}|\xi^2_{\ell}|^{2<\ell<m-3}&=&p_{m-5}|\xi^2_4|+p_{m-7}|\xi^2_6|+\cdots+p_{6}|\xi^2_{m-7}|+p_{4}|\xi^2_{m-5}|\notag\\
       &\leq&(C_{m-5}C_4+C_{m-7}C_6+\cdots+C_6C_{m-7}+C_4C_{m-5})\,\beta^2J^{m-1}\notag\\
       &=&(C_{m-1-\ell}C_{\ell}^{(\ell\,\,even)}-2C_2C_{m-3}-2C_0C_{m-1})\,\beta^2J^{m-1}\notag\\
       &\leq&(E(4)C_{m-1-\ell}C_{\ell}-2C_2C_{m-3}-2C_0C_{m-1})\,\beta^2J^{m-1}\notag\\
       &=&(E(4)C_m-2C_2C_{m-3}-2C_0C_{m-1})\,\beta^2J^{m-1}\notag\\
       &=&\left(\beta\frac{E(4)C_m-2C_2C_{m-3}-2C_0C_{m-1}}{C_{m-1}}\right)\,\beta J^{m-1}C_{m-1}\notag\\
       &\leq&\beta(E(4)4-1/4-2)\,\beta J^{m-1}C_{m-1},\notag\\
       &\leq&(0.7976\beta)\,\beta J^{m-1}C_{m-1},\label{xithree}
\end{eqnarray}
Thus, adding \eqref{xione}, \eqref{xitwo} and \eqref{xithree}, we set
\begin{eqnarray}
S^*=4J\{\thetapbar Q^*+\thetap(E(4)\beta J^{-2}(1/3)+E^2(4)\beta^2J^{-3}(1/3)
+E(4)\beta J^{-3}(5/42))&&\notag\\
+\thetapbar(E(4)\beta J^{-2}(1/3)+E(4)\thetapbar\beta J^{-3}(5/42)+\thetap J^{-3}(1/21))\}+\notag\\
+\thetap\{(\,|\xi^2_0|R^*+|\xi^2_2|J^{-2}(1/7)+p_2J^{-2}(1/7)\,)+(0.7976\beta)\}.\label{sstar}
\end{eqnarray}
Thus, taking $J_2$ such that $Q^*(J_2),\,R^*(J_2),\,S^*(J_2)\leq 1$ and $J=\text{max}\{J_1,\,J_2\}$, we have shown that the induction hypothesis \eqref{catalanhypothesis} implies
\begin{equation}\label{finalcatalan}
\bar p_m,\,p_m,\,|\xi^2_m|\leq \beta C_mJ^m,
\end{equation}
so that in fact \eqref{finalcatalan} holds for every integer $m$.
\end{proof}

\subsection{Proof of Theorem \ref{thm:solution}: Solution of the Eigenvalue Problem}\label{solutionproof}
\begin{proof}
The weak form of the master system is
\begin{equation}\label{weakmaster}
  \int_Q \left[ \eps_\eta^{-1}(\nabla+i\eta\kk)h_\eta(y)\cdot(\nabla-i\eta\kk)\bar v(y) - \eta^2 \xi^2_\eta h_\eta(y) \bar v(y) \right] = 0
  \quad \text{for all } v\in H^1_\text{\tiny per}(\QQ).
\end{equation}
Using that $\eps_\eta^{-1} = \eta^2\xi^2_\eta/(\eta^2\xi^2_\eta-1)$ in $P$ and $\eps_\eta=1$ in $\Pbar$, and multiplying by $(\eta^2\xi^2_\eta-1)$, gives the equivalent system
\begin{equation*}
  a_\eta(h,\xi^2;v) = 0 \quad \text{for all } v\in H^1_\text{\tiny per}(\QQ),
\end{equation*}
in which
\begin{multline*}
  a_\eta(h,\xi^2;v) = - \int_\Pbar (\nabla+i\eta\kk)h\cdot(\nabla-i\eta\kk)\bar v \; + \\
  + \int_\QQ \left[ \eta^2\xi^2 (\nabla+i\eta\kk)h\cdot(\nabla-i\eta\kk)\bar v - (\eta^2\xi^2-1) \eta^2\xi^2h\bar v \right].
\end{multline*}
This form can be expanded in powers of $\eta$,
\begin{equation}\label{aexpansion}
  a_\eta(h,\xi^2;v) = a_0(h;v) - i\eta a_1(h;v) - \eta^2 a_2(h,\xi^2;v) + i\eta^3 a_3(h,\xi^2;v) + \eta^4 a_4(h,\xi^2;v) 
\end{equation}
in which the $a_m$ are real forms
\begin{eqnarray*}
  && a_0(h;v) = -\int_\Pbar \nabla h\cdot\nabla\bar v, \\
  && a_1(h;v) = \int_\Pbar (h\kk\cdot\nabla\bar v - \nabla h\cdot\kk\bar v), \\
  && a_2(h,\xi^2;v) = \int_\Pbar h\bar v - \xi^2\int_\QQ(\nabla h\cdot\nabla\bar v + h\bar v), \\
  && a_3(h,\xi^2;v) = \xi^2 \int_\QQ (h\kk\cdot\nabla\bar v - \nabla h\cdot\kk\bar v), \\
  && a_4(h,\xi^2;v) = \xi^2\int_\QQ (1-\xi^2)h\bar v.
\end{eqnarray*}
Define the partial sums
\begin{equation*}
  \xi^{2,N}_\eta = \sum_{m=0}^N \eta^m \xi^2_m\,\,\,\,\,\,\mbox{and}\,\,\,\,\,\,h^N_\eta = \sum_{m=0}^N \eta^m h_m.
\end{equation*}
For $\eta<R$,
the sequence $\{\xi^{2,N}_\eta\}$ converges to a number $\xi^2_\eta$ and the sequence $\{h^N_\eta\}$ converges in $H^1_*$ to a function $h_\eta$; thus

\begin{equation*}
  a_j(h^N_\eta,\xi^{2,N}_\eta;v) \to a_j(h_\eta,\xi^2_\eta;v) \quad \text{for all } v\in H^1_\text{\tiny per}(\QQ)  \text{ and } i=0,\dots,4.
\end{equation*}
Therefore, $a_j(h_\eta,\xi^2_\eta;v)$, $j=1,\dots,4$, has a convergent series representation in powers of $\eta$, in which the $m^\text{th}$ coefficient is related to the coefficients $\xi_\ell$ and $h_\ell$ by
\begin{eqnarray*}
  &(j=0)& \int_\Pbar -\nabla h_m\cdot\nabla\bar v, \\
  &(j=1)& \int_\Pbar \left( h_m\kk\cdot\nabla\bar v - \nabla h_m\cdot\kk\bar v \right), \\
  &(j=2)& \int_\Pbar h_m\bar v - \int_\QQ \left( \nabla(\xi^2_\ell h_{m-\ell})\cdot\nabla\bar v + (\xi^2_\ell h_{m-\ell})\bar v \right), \\
  &(j=3)& \int_\QQ \left( (\xi^2_\ell h_{m-\ell})\kk\cdot\nabla\bar v - \nabla(\xi^2_\ell h_{m-\ell})\cdot\kk\bar v \right), \\
  &(j=4)& \int_\QQ \left( \xi^2_\ell h_{m-\ell} - \xi^2_j\xi^2_{\ell-j} h_{m-\ell} \right)\bar v.
\end{eqnarray*}
From these, one obtains the $m^\text{th}$ coefficient of $a_\eta(h_\eta,\xi^2_\eta;v)$ (see \ref{aexpansion}), which, by means of the relations $h_m= h_0i^m\psi_m$, $\xi^2_\ell h_{m-\ell} = h_0i^m \sigma'_m$ and $\xi^2_j\xi^2_{\ell-j} h_{m-\ell} = h_0 i^m \sigma''_m$ is seen to be equal to the $-i^mh_0$ times the right-hand side of equation \eqref{weakform}.  All these coefficients are therefore equal to zero, and we conclude that $a_\eta(h_\eta,\xi^2_\eta;v) = 0$. This proves that the function $h_\eta$, together with the frequency $\sqrt{\xi_\eta}$ solve the weak form \eqref{weakmaster} of the master system.
\end{proof}
\section{Effective Properties, Error Bounds and the Dispersion Relation}\label{sec:errorbounds}

In this section we start by identifying an effective property directly from the dispersion relation. We then discuss the relation between effective properties and quasistatic properties. Next we provide explicit error bounds for finite-term approximations to the first branch of the dispersion relation for nonzero values of $\eta$. The error bounds show that numerical computation of the first two terms of the power series delivers an accurate and inexpensive numerical method for calculating dispersion relations for sub-wavelength plasmonic crystals. 

\subsection{The Effective Index of Refraction - Quasistatic Properties and Homogenization}
\label{effproperties}

The identification of an effective index of refraction valid for $\eta>0$  follows directly from the dispersion relation given by the series for $\omega^2_\eta$. Indeed the effective refractive index $n^2_{\text{eff}}$ is defined by expressing  the dispersion relation as 
\begin{eqnarray}
\omega_\eta^2=\frac{c^2 k^2}{n^2_{\text{eff}}}
\label{effeispersdyn}
\end{eqnarray}
and it then follows from the expansion for $\omega^2_\eta$ that 
the effective refractive index has the convergent power series expansion
\begin{eqnarray}
n_{\text{eff}}^{-2}=n^{-2}_{\text{qs}}+\sum_{m=1}^\infty\eta^{2m}\xi^2_{2m}.\label{expandindex}
\end{eqnarray}

We now discuss the relationship between the effective index of refraction and 
the quasistatic effective properties seen in the $d\rightarrow 0$ limit with $k$ fixed. The effective refraction index $n_{\text{eff}}$ can be rewritten in the equivalent form by the equation
$n^2_{\text{eff}}=1/\xi^2_{\eta}$.
By setting $v=h_{\eta}$ in the weak form of the master system \eqref{weakmaster}, it is easily seen that $\xi^2_{\eta}>0$ for all $\eta$ within the radius of convergence, so that $n^2_{\text{eff}}>0$ for those values of $\eta$. Following Pendry \textit{et al.} (1999), see also Kohn \& Shipman (2008), we  define the effective permeability by 

\begin{equation}
\mueff=\frac{(B_3)_{\text{eff}}}{(H_3)_{\text{eff}}},
\label{mueff}
\end{equation}
and we then define $\epseff$ through the equation
\begin{equation}\label{epseff}
n^2_{\text{eff}}=\epseff\mueff.
\end{equation}
The quasi-static effective properties are recovered by passing to the limits
\begin{eqnarray*}
n^2_\text{qs}=\lim_{\eta\rightarrow 0}n^2_{\text{eff}},\,\,\,\,\mu_\text{qs}=\lim_{\eta\rightarrow 0}\mu_\text{eff},\,\,\,\,
\eps_\text{qs}=\lim_{\eta\rightarrow 0}\epseff.
\end{eqnarray*}
A simple computation shows that $\mu_\text{qs}=\integraloverq{\psi_0}>0$ (see appendix). Hence, we have that $\mueff>0$ for $\eta$ in a neighborhood of the origin, so that $\epseff>0$ for these values of $\eta$, since $n^2_{\text{eff}}>0$ for all $\eta$ in the radius of convergence. Thus, \emph{one has a solid basis on which to assert that plasmonic crystals function as materials of positive index of refraction in which both the effective permittivity and permeability are positive.}

For circular inclusions we have used the program COMSOL to compute that $\integraloverq{\psi_0}\approx 0.98$, so that only a mild effective magnetic permeability arises.

Having established that $h_{\eta}(\yy)e^{i\eta\kk\cdot \yy}$ is the solution to the unit cell problem, we can undo the change of variable $\yy=k\xx/\eta$ to see that the function
\begin{equation}
\hat h_\eta\left(\frac{k\xx}{\eta}\right)e^{i(k\kk\cdot\xx-\omega_{\eta}t)},
\label{crystalsoln}
\end{equation}
where $\hat h_{\eta}$ is the $\QQ$-periodic extension of $h_{\eta}$ to all of $\RR^2$, is a solution of
\begin{equation}\label{macrowave}
\nabla_\xx\cdot(\eps^{-1}_{\eta}\nabla_\xx \hat h_{\eta})=\frac{1}{c^2}\partial_{tt}\hat h_{\eta},\,\,\,\,\xx\in\RR^2,
\end{equation}
for every $\eta$ in the radius of convergence.

We investegate the quasistatic limit directly using the power series \eqref{crystalsoln}. Here we wish to describe the average field as $d\rightarrow 0$.
To do this we introduce the three-dimensional period cell for the crystal $[0,d]^3$. The base of the cell in the $x_1x_2$ plane is denoted by $Q_d=[0,d]^2$ and is the period of the crystal in the plane transverse to the rods. We apply the definition of $B_\text{eff}$ and $H_\text{eff}$ given in Pendry \textit{et al.} (1999) which in our context is
\begin{eqnarray}
(B_3)_{\text{eff}}=\frac{1}{d^2}\int_{Q_d}\hat h_\eta\left(\frac{k\xx}{\eta}\right)e^{i(k\kk\cdot\xx-\omega_{\eta}t)}dx_1 dx_2
\label{avgbfield}
\end{eqnarray}
and
\begin{eqnarray}
(H_3)_{\text{eff}}=\frac{1}{d}\int_{(0,0,0)}^{(0,0,d)}\hat{h}_\eta\left(\frac{k\xx}{\eta}\right)e^{i(k\kk\cdot\xx-\omega_{\eta}t)}dx_3.
\label{avghfield}
\end{eqnarray}
Taking limits for $k$ fixed and $d\rightarrow 0$ in \eqref{crystalsoln} gives
\begin{eqnarray*}
\lim_{d\rightarrow 0}(B_3)_{\text{eff}}=\langle\psi_0\rangle_{\scriptscriptstyle{\QQ}}\bar{h}_0 e^{i(k\kk\cdot\xx-\omega_\text{qs}t)} &\hbox{  and  }&\lim_{d\rightarrow 0}(H_3)_{\text{eff}}=\bar{h}_0 e^{i(k\kk\cdot\xx-\omega_\text{qs}t)},
\label{avgbandh}
\end{eqnarray*}
in which $\omega^2_{\text{qs}}=\frac{c^2k^2}{n^2_{\text{qs}}}$.  These are the same average fields that would be seen in a quasistatic magnetically active effective medium with index of refraction $n_\text{qs}$ and $\mu_{\text{qs}}$ that supports the plane waves
\begin{eqnarray*}
(B_3)_{\text{qs}}=\mu_{qs}\bar{h}_0e^{i(k\kk\cdot\xx-\omega_{\text{qs}}t)}  &\hbox{  and  }(H_3)_{\text{qs}}=\bar{h}_0e^{i(k\kk\cdot\xx-\omega_{\text{qs}}t)},
\label{avgbandb}
\end{eqnarray*}
where $\mu_{\text{qs}}=\langle \psi_0\rangle_\QQ$. It is evident that these fields are solutions of the homogenized equation
$\frac{n^2_{\text{qs}}}{c^2}\partial_{tt}u=\Delta u$. This quasistatic interpretation provides further motivation for the  definition  of  $n_{\text{eff}}$ for nonzero $\eta$ given by \ref{expandindex}.

Now we apply the definition of effective permeability $\mu_{\text{eff}}$ given in Pendry \textit{et al.} (1999),   together with $n_{\text{eff}}$ to define an effective permeability $\epsilon_{\text{eff}}$ for $\eta>0$.  The relationships between the effective properties and quasistatic effective properties are used to show that plasmonic crystals function as meta-materials of positive index of refraction in which both the effective permittivity and permeability are positive for $\eta>0$.
\subsection{Absolute Error Bounds}
The Catalan bound provides simple estimates on the size of the tails for the series $\xi^2_{\eta}$ and~$h_{\eta}$,
\begin{eqnarray*}
E_{m_0,\xi}=\sum_{m=m_0+1}^{\infty}\xi^2_{2m}\eta^{2m}\,\,\,\,\,\,\mbox{and}\,\,\,\,\,\,E_{m_0,h}&=&\sum_{m=m_0+1}^{\infty}h_m\eta^m.
\end{eqnarray*}
We have established convergence for $\eta\leq 1/4J$, so that we may write $\eta=\alpha/4J$, $0\leq\alpha\leq 1$. Then, using that $|\xi^2_{2m}|\leq \beta C_{2m}J^{2m}$, $C_{2m}\leq 4^{2m}$ and $4J\eta=\alpha$, we have
\begin{eqnarray}
		|E_{m_0,\xi}|&=&\left|\sum_{m=m_0+1}^{\infty}\xi^2_{2m}\eta^{2m}\right|\notag\\
                   &\leq&\beta\sum_{m=m_0+1}^{\infty}C_{2m}J^{2m}\eta^{2m}\notag\\
									 &\leq&\beta\sum_{m=m_0+1}^{\infty}(4J\eta)^{2m}\notag\\
									 &\leq&\beta\sum_{m=m_0+1}^{\infty}\alpha^{2m}\notag\\
									 &=   &\beta\frac{\alpha^{2m_0+2}}{1-\alpha^2}.\label{xierror}
\end{eqnarray}
Similarly, for $h_{\eta}$ we have that
\begin{eqnarray}
     \Vert E_{m_0,h}\Vert_{H^1(\QQ)}
																		&\leq&2\beta|h_0|\frac{\alpha^{m_0+1}}{1-\alpha}.\label{herror}
\end{eqnarray}
\subsection{Relative Error Bounds}\label{subsec:errorbounds}
In this section, we use the absolute error bound \eqref{herror} with $m_0=1$ to obtain a relative error bound for the particular case of a circular inclusion of radius $r=0.45$ Shvets \& Urzhumov (2004). The first term approximation to $h_{\eta}$ is
\begin{equation}
h_{\eta}=h_0\psi_0+ih_0\psi_1\eta+E_{1,h}\label{happrox}.
\end{equation}
For a circular inclusion of radius $r=0.45$, we have
$
J\leq85,\,\beta\leq0.79,\,\Vert\psi_0\Vert=0.97\,\text{and}\,\Vert\psi_1\Vert=0.02,
$
where $\Vert\cdot\Vert=\Vert\cdot\Vert_{H^1(\QQ)}$. Thus, using bound \eqref{herror}, the relative error $R_{1,h}$ is bounded by
\begin{multline*}
|R_{1,h}|=\frac{|E_{1,h}|}{\Vert h_0+ih_0\psi_1\eta\Vert}
         \leq \frac{1.58\frac{\alpha^2}{1-\alpha}}{\Vert h_0\Vert-\Vert h_0\psi_1\Vert |\eta|}
         \leq \frac{1.58\frac{\alpha^2}{1-\alpha}}{0.97-0.02\frac{\alpha}{340}},
\end{multline*} 
so that for $\alpha\leq 0.2$ the relative error is less than $8.2\%$. The graphs of $\psi_0$ and $\psi_1$ can be found in the Appendix. The first term approximation to $\xi^2_{\eta}$ is
\begin{equation}
\xi^2_{\eta}=\xi^2_0+\xi^2_2\eta^2+E_{1,\xi}.\label{xiapprox}
\end{equation}
In the Appendix we indicate how the tensors $\xi^2_m$ may be computed. For an inclusion of radius $r=0.45$, we have $\xi^2_0=0.36$ and $\xi^2_2=-0.14$. Thus, using bound \eqref{xierror}, the relative error $R_{1,\xi}$ is bounded by
\begin{eqnarray*}
|R_{1,\xi}|=\frac{|E_{1,\xi}|}{|\xi^2_0\psi_0+\xi^2_2\eta^2|}
\leq\frac{\beta\frac{\alpha^4}{1-\alpha^2}}{|\xi^2_0+\xi^2_2\eta^2|}
\leq
\frac{0.79\frac{\alpha^4}{1-\alpha^2}}{|0.36-0.14\frac{\alpha^2}{340^2}|},
\end{eqnarray*}
so that for $\alpha\leq0.3$ the relative error is less than $2\%$.
\begin{figure}
\begin{minipage}[b]{0.5\linewidth}
\centering
\includegraphics[scale=0.25]{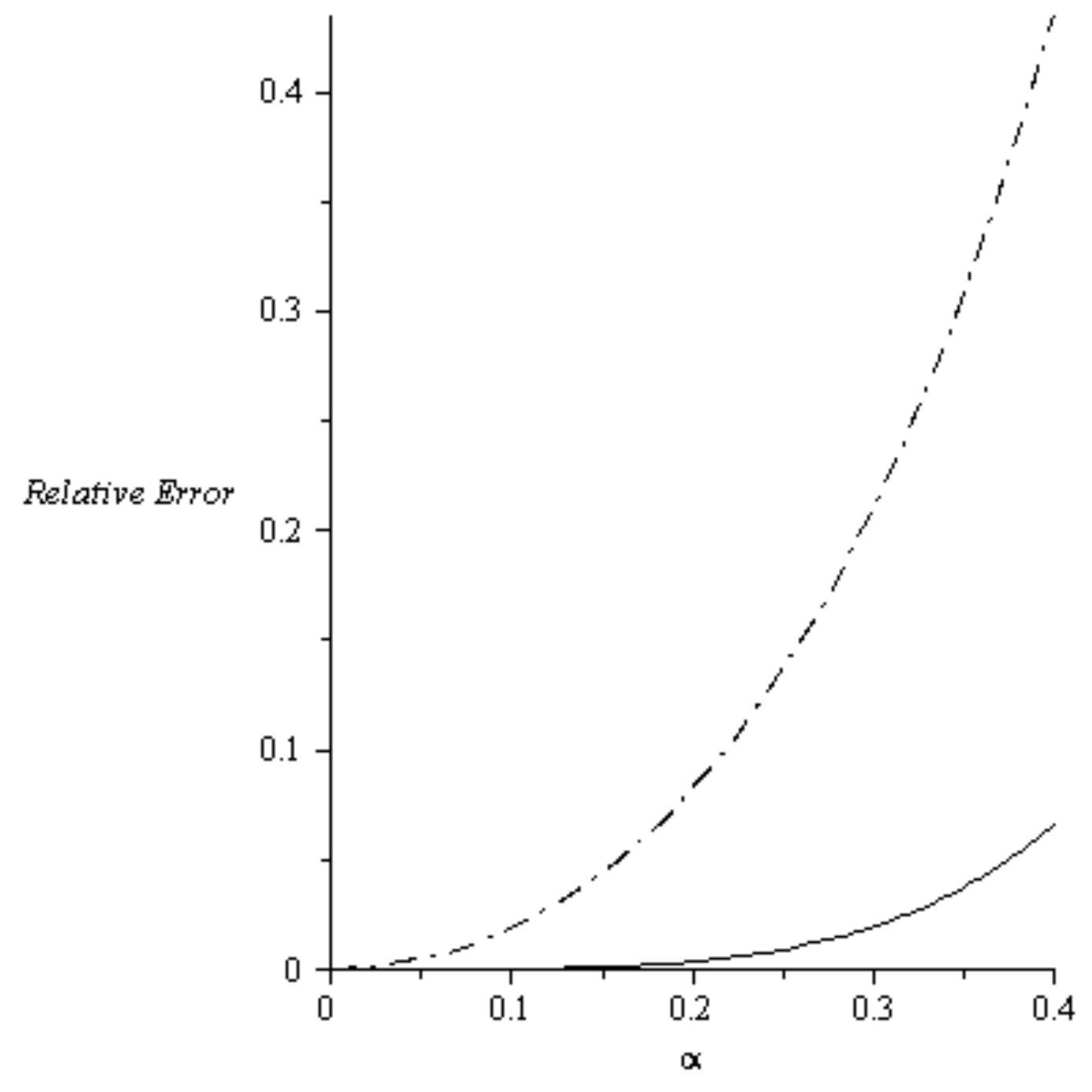}
\caption{Solid curve is $R_{1,\xi}$ and dotdash curve is $R_{1,h}$.}
\label{fig:relativeerror}
\end{minipage}
\begin{minipage}[b]{0.5\linewidth}
\centering
\includegraphics[scale=0.25]{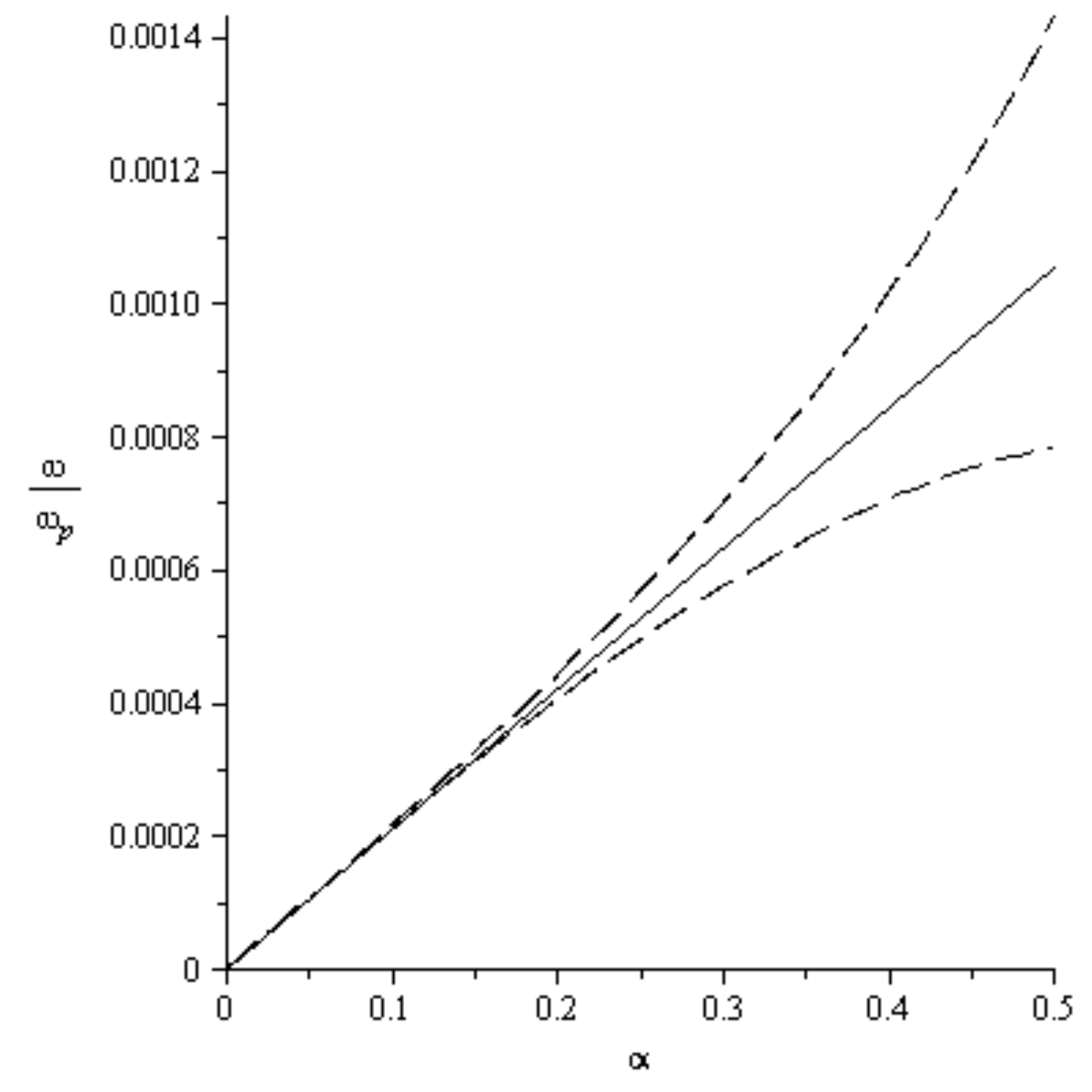}
\caption{Graph of the first branch of the dispersion relation. Dashed curves represent error bars.}
\label{fig:dispersionbranch}
\end{minipage}
\end{figure}
\begin{acknowledgements}
R. Lipton was supported by NSF grant DMS-0807265 and AFOSR grant FA9550-05-0008, and S. Shipman was supported by NSF grant DMS-0807325.  The Ph.D. work of S. Fortes was also supported by these grants.  This work was inspired by the IMA ``Hot Topics" Workshop on Negative Index Materials in October, 2006.
\end{acknowledgements}
\appendix{Explicit Expressions for Tensors}
The tensors $\xi^2_0$ and $\xi^2_2$ were calculated using the weak form \eqref{weakform}, as follows. Setting $m=2$ and $v\equiv 1$ in \eqref{weakform} and solving for $\xi^2_0$ gives
\begin{equation}
\xi^2_0=\frac{\integraloverpbar{\kk\cdot\!\!\nabla\psi_1+\psi_0}}{\integraloverq{\psi_0}}.\label{xizerosquare}
\end{equation}
Setting $m=4$ and $v\equiv 1$ and solving for $\xi^2_2$ gives
\begin{equation}
\xi^2_2=\frac{-\xi^2_0\integraloverq{\kk\cdot\!\!\nabla\psi_1}+\xi^2_0\integraloverp{\psi_2}+\xi^2_0\xi^2_0\integraloverq{\psi_0}+\xi^2_0\integraloverq{\psi_0}-\vectorintegraloverpbar{\kk}{\nabla\psi_3}}{\integraloverq{\psi_0}}.\label{xitwosquare}
\end{equation}
All integrals appearing in \eqref{xizerosquare} and \eqref{xitwosquare} were then computed using the program COMSOL.
\appendix{Bounds on $\bar p_0$, $p_0$ and $|\xi^2_0|$}
We have that $\bar p_0=\theta_{\Pbar},\,p_0\leq\thetap,\text{ and }\xizero\leq 1$. These bounds are obtained as follows: we have $\psi_0\equiv 1$ in $\bar P$, so that $\bar p_0=\theta_{\Pbar}$. Using the BVP for $\psi_0$ in $P$ one can prove that $0\leq\psi_0(\yy)\leq1$, $\forall \yy\in P$, and that $p_0^2=\integraloverp{\psi_0}$. These two facts together give the estimate $p_0\leq\thetap$. Setting $v=\psi_1$ in the weak form for $\psi_1$, we get that $\vectorintegraloverpbar{\kk}{\nabla\psi_1}=-\vectorintegraloverpbar{\nabla\psi_1}{\nabla\psi_1}<0$. Using this in expression \eqref{xizerosquare} gives $\xi^2_0<1$. Since $\theta_P,\,\theta_{\bar P}\leq 1$, these three estimates allow us to take $\beta\leq1$.

\appendix{Computing the Constant A for Circular Inclusions}

\newcommand{\II}{\mathrm{II}}
\newcommand{\KK}{\mathrm{KK}}
\newcommand{\IK}{\mathrm{IK}}
\newcommand{\KI}{\mathrm{KI}}
\newcommand{\JJ}{\mathrm{JJ}}

Given a function $\psi\in H^1_{\text{per}}(\bar P)$, let $u\in H^1(P)$ satisfy
\begin{equation*}
  \renewcommand{\arraystretch}{1}
\left.
  \begin{array}{ll}
    \nabla^2 u - u = 0 & \text{in }P,\\
    u=\psi & \text{on } \partial P.
  \end{array}
\right.
\end{equation*}
We seek to compute a number $A$ such that $\|u\|^2_{H^1(P)} \leq A\|\psi \|^2_{H^1(\bar P)}$ for all $\psi$.  Following Bruno (1991), we will calculate a value of $A$ for circular inclusions $P$ of radius $r_0<0.5$ by restricting $\psi$ to the annulus between $P$ and the circle of unit radius.  It suffices to consider real-valued functions $\psi$ that minimize the $H^1$ norm in the annulus, that is $\nabla^2\psi-\psi=0,\quad r_0<r<0.5$. A function of this type is given generally by the real part of an expansion $\psi(r,\theta) = \sum_{n=0}^\infty \left(c_nI_n(r) + d_nK_n(r)\right) e^{in\theta}$, in which $c_n$ and $d_n$ are complex numbers and $I_n$ and $K_n$ are the ``modified" Bessel functions.  The continuous continuation of $\psi$ into the disk with $\nabla^2 u - u = 0$ is given by the real part of $u(r,\theta) = \sum_{n=0}^\infty f_nI_n(r) e^{in\theta}$ under the relations
\begin{equation}\label{matchingcoeff}
  f_n = c_n + d_n\frac{K_n(r_0)}{K_n(r_0)}.
\end{equation}
One computes that $\|\Re\psi\|^2=\frac{1}{2}\|\psi\|^2\quad\text{and}\quad \|\Re u\|^2=\frac{1}{2}\|u\|^2$, so we may work with the complex functions rather than their real parts. The Helmholtz equation in $P$ and integration by parts yield
\begin{equation*}
  \|u\|_{H^1(P)}^2 = \int_P \left(|\nabla u|^2 + |u|^2\right)dA =
    \int_{\partial P} \bar u \partial_n u,
\end{equation*}
and this provides the representation $\|u\|^2_{H^1} = 2\pi r_1 \sum_{n=0}^\infty \bar f_n I_n(r_1) f_nI'_n(r_1)$. The analogous representation in the annulus is
\begin{multline*}
  \|\psi\|^2_{H^1} = \int_0^{2\pi} (\partial_r\psi(1,\theta)\overline{\psi(1,\theta)}\,d\theta
                 - \int_0^{2\pi} r_0(\partial_r\psi(r_0,\theta)\overline{\psi(r_0,\theta)}\,d\theta \\
   =\; 2\pi \sum_{n=0}^\infty \left(\bar c_nI_n(0.5)+\bar d_n K_n(0.5)\right)\left(c_n I'_n(0.5)+d_nK'_n(0.5)\right) + \\
      - 2\pi r_0\sum_{n=0}^\infty \left(\bar c_nI_n(r_0)+\bar d_n K_n(r_0)\right)\left(c_n I'_n(r_0)+d_nK'_n(r_0)\right)
\end{multline*}

We seek a positive number $A$ such that, for all choices of complex numbers $c_n$ and $d_n$ we have $ 0 \leq A\|\psi\|^2 - \|u\|^2$. The right-hand-side of this inequality is a quadratic form in all of the coefficients $(c_n,d_n)$ that depends on $A$,
\begin{equation*}
 A\|\psi\|^2 - \|u\|^2  = m_n^{11} c_n\bar c_n + m_n^{12} c_n\bar d_n + m_n^{21} d_n \bar c_n + m_n^{22} d_n\bar d_n,
\end{equation*}
in which the $m_n^{ij}$ depend on $A$ and are conveniently expressed in terms of the functions
\begin{eqnarray*}
  \II_n(r)=I_n(r)I'_n(r),&\KK_n(r)=K_n(r)K'_n(r),&\IK_n(r)=I_n(r)K'_n(r)\\
  \KI_n(r)= K_n(r)I'_n(r),&\JJ_n(r)=\frac{K_n(r)^2}{I_n(r)}I'_n(r).
\end{eqnarray*}
\begin{eqnarray*}
  m_n^{11} &=& -r_0\II_n(r_0) - Ar_0\II_n(r_0) + A\,\II_n(0.5), \\
  m_n^{22} &=& -r_0\JJ_n(r_0) - Ar_0\KK_n(r_0) + A\,\KK(0.5), \\
  m_n^{12} &=& -r_0\KI(r_0) - Ar_0\KI_n(r_0) + A\,\KI_n(0.5), \\
  m_n^{21} &=& -r_0\KI(r_0) - Ar_0\IK_n(r_0) + A\,\IK_n(0.5).  
\end{eqnarray*}
The form $m_n^{ij}$ is Hermitian, as one can show that $m_n^{12}=m_n^{21}$ by using the fact that $r\text{Wron}[I_n,K_n]$ is constant.

We must find $A>0$ such that $m^{11}_n\geq0$ and $m_n^{11}m_n^{22}-m_n^{12}m_n^{21}\geq0$ for all $n=0,1,2,\dots$.  These quantities are equal to
\begin{eqnarray*}
  m_n^{11} &=& \beta(r_0,0.5) A - \alpha(r_0), \\
  m_n^{11}m_n^{22}-m_n^{12}m_n^{21} &=& \epsilon(r_0,0.5) A^2 - \delta(r_0,0.5) A,
\end{eqnarray*}
in which
\begin{eqnarray*}
  \alpha_n(r) &=& r\II_n(r), \\
  \beta_n(r,s) &=& \II_n(s) - r\II_n(r), \\
  \delta_n(r,s) &=& rs\left[
         \II_n(r)\KK_n(s) + \II_n(s)\JJ_n(r) - \KI_n(r)\IK_n(s) - \KI_n(s)\KI_n(r) \right] + \\
                       && +\, r^2\left[
         -\II_n(r)\KK_n(r) - \II_n(r)\JJ(r) + \KI_n(r)\IK_n(r) +\KI_n(r)\KI_n(r) \right], \\
  \epsilon_n(r,s) &=& rs\left[ 
         -\II_n(r)\KK_n(s) - \II_n(s)\KK_n(r) + \KI_n(r)\IK_n(s) + \KI_n(s)\IK_n(r) \right].
\end{eqnarray*}
The numbers $\alpha_n(r)$ and $\beta_n(r,s)$ for $r<s$ are positive; the latter because
$
  (rI_nI'_n)' = \frac{1}{r}(r^2+n^2)I^2_n + r{I'_n}^2 > 0.
$
One can show that $\delta$ and $\epsilon$ are positive.
Thus, it is sufficient to find $A>0$ such that, for all $n=0,1,2,\dots$,
\begin{equation*}
  A \geq \max\left\{\frac{\alpha_n(r_0)}{\beta_n(r_0,0.5)}, \frac{\delta_n(r_0,0.5)}{\epsilon_n(r_0,0.5)} \right\}.
\end{equation*}
Table \ref{Jtable} shows computed values of $A$ for various values of $r_0$.
\appendix{Graphs of $\psi_m$ and Table of $A$, $\Omegap$ and $J$ for Circular Inclusions}
\begin{figure}[ht]
\begin{minipage}[b]{0.5\linewidth}
\centering
\includegraphics[scale=0.15]{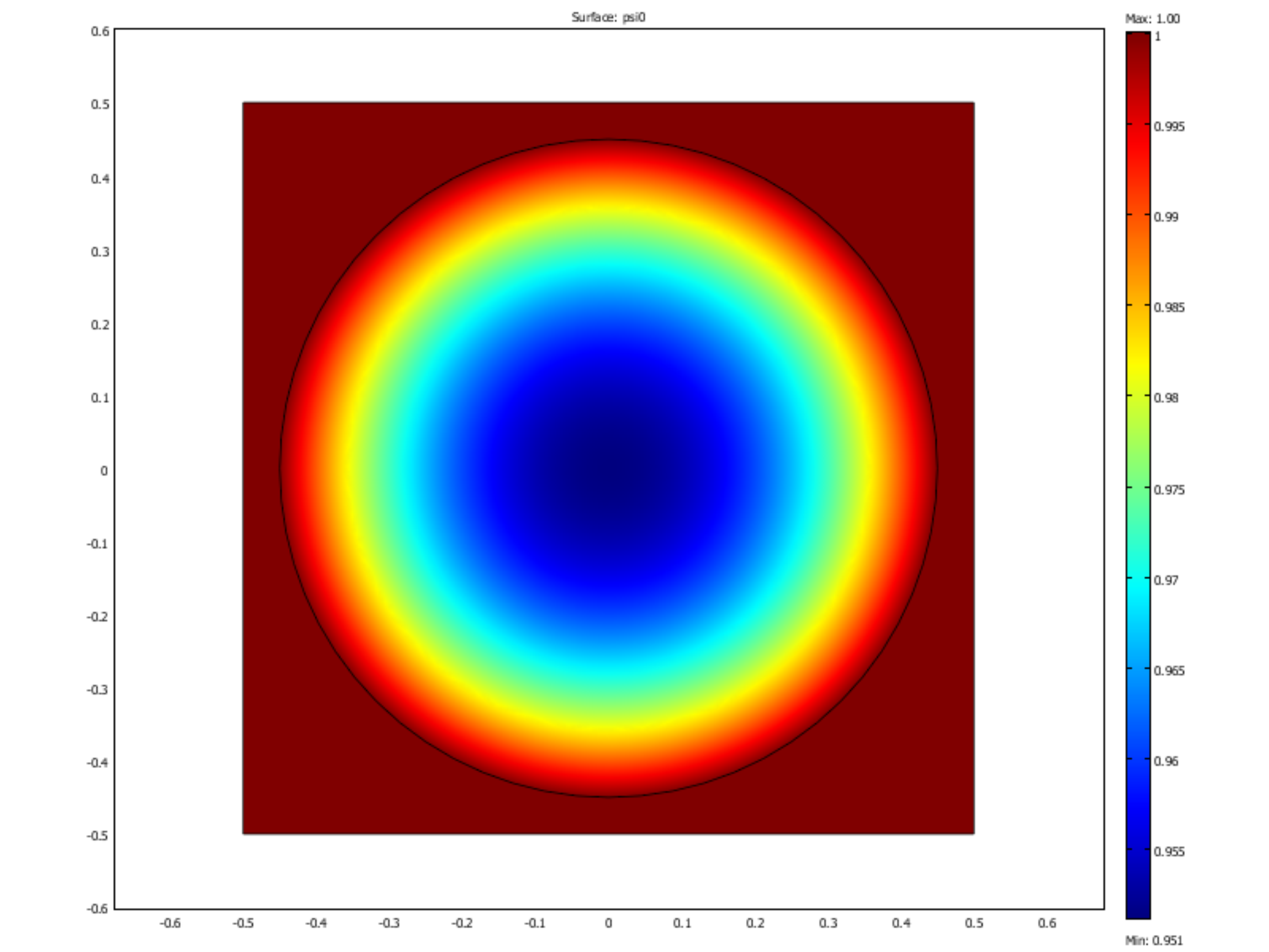}
\caption{Graph of $\psi_0$.  This function is symmetric about the origin.}
\label{fig:figure1}
\end{minipage}
\begin{minipage}[b]{0.5\linewidth}
\centering
\includegraphics[scale=0.15]{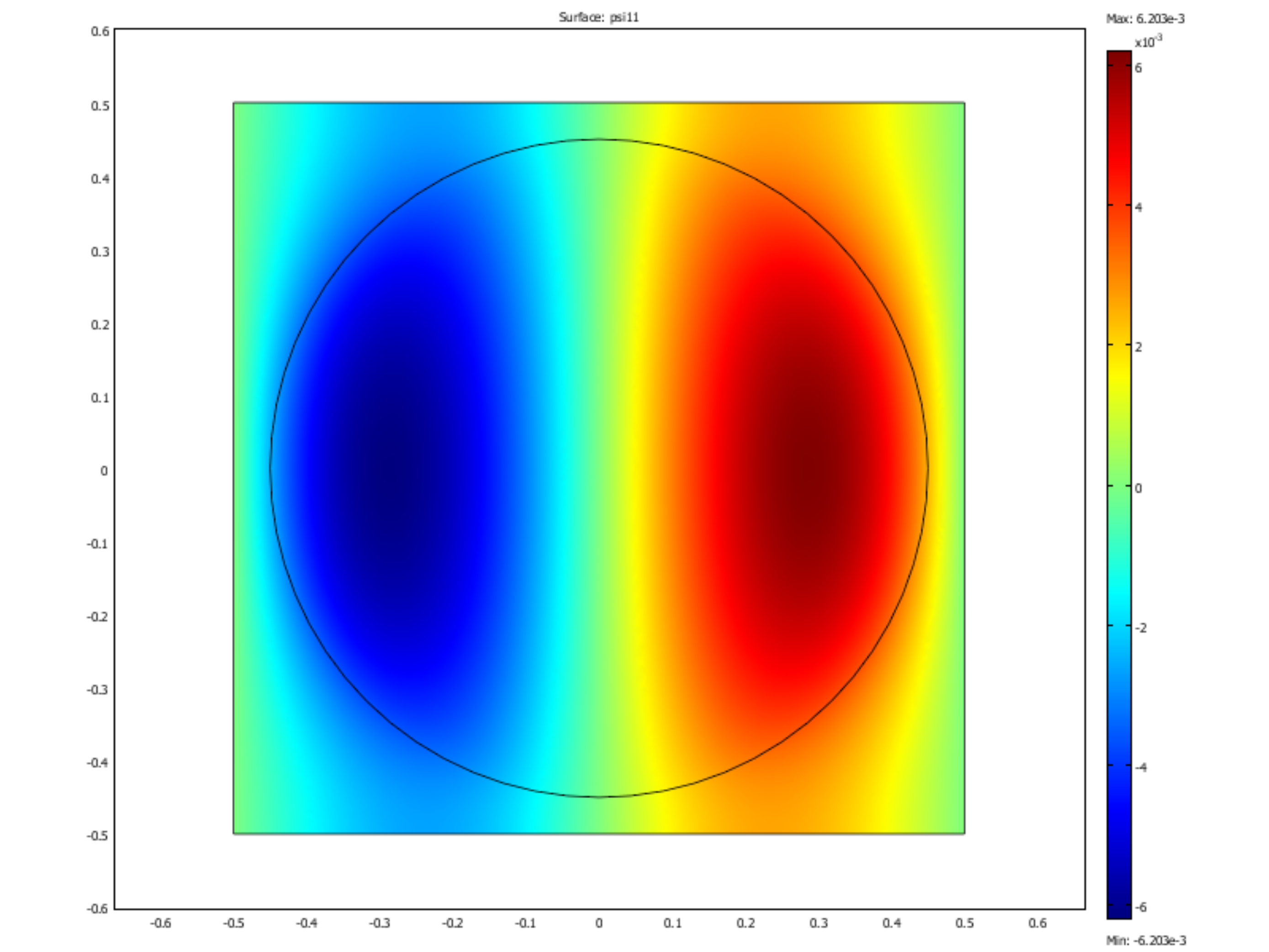}
\caption{Graph of $\psi_1$ when $\kk=(1,0)$.  This function is antisymmetric about the origin.}
\label{fig:figure2}
\end{minipage}
\end{figure}
%
%
\begin{table}\label{Jtable}
\centering 
\begin{tabular}{|c|c|c|c|c|c|}
\hline
      $r$       &0.1    &0.2  &0.3    &0.4    &0.45    \\ \hline
      $A$       & 1.058 &1.293&1.907  &3.956  & 4.840 \\ \hline
      $\Omegap$ & 1.0     & 1.0   & 1.0     & 1.0     &1.0 \\ \hline
      $J$       & 15    &17   &22     &29     &85 \\ \hline
\end{tabular}
\caption{Values of $A$, $\Omegap$ and $J$ for circular inclusions of radii $rd$.}
\end{table}

\newpage

\medskip
\noindent
\center {{\large \bfseries Figure Captions.}} \
\begin{itemize}
\item[1] Unit cell with plasmonic inclusion.

\item[2] Solid curve is $R_{1,\xi}$ and dotdash curve is $R_{1,h}$.

\item[3] Graph of $\psi_0$.  This function is symmetric about the origin.

\item[4] Graph of $\psi_1$ when $\kk=(1,0)$.  This function is antisymmetric about the origin.
\end{itemize}

\medskip
\noindent
\center {{\large \bfseries Short Title for Heading: Subwavelength Plasmonic Crystals}} \

\end{document}